\documentclass[12pt,a4paper]{article}

\usepackage{pdflscape}
\usepackage{amsmath}
\usepackage{amssymb}
\usepackage{latexsym}
\usepackage{srcltx}
\usepackage{graphics}
\usepackage{color}
\usepackage{epsfig}
\usepackage{color}
\usepackage{url}

\newcommand{\re}{{\mathbb R}}
\newcommand{\n}{{\mathbb N}}

\newcommand{\cA}{{\mathcal{A}}}

\newcommand{\cV}{{\mathcal{V}}}

\newcommand{\cT}{{\mathcal{T}}}

\newcommand{\cS}{{\mathcal{S}}}
\newcommand{\cR}{{\mathcal{R}}}

\newcommand{\cP}{{\mathcal{P}}}

\newcommand{\cU}{{\mathcal{U}}}

\newcommand{\bx}{{\boldsymbol{x}}}
\newcommand{\by}{{\boldsymbol{y}}}

\newcommand{\ba}{{\boldsymbol{a}}}

\newcommand{\bc}{{\boldsymbol{c}}}
\newcommand{\bv}{{\boldsymbol{v}}}
\newcommand{\bu}{{\boldsymbol{u}}}

\newcommand{\bn}{{\boldsymbol{n}}}
\newcommand{\bh}{{\boldsymbol{h}}}

\newcommand{\vardot}{\mathord{\,\cdot\,}}

\newtheorem{theorem}{Theorem}
\newtheorem{prop}{Proposition}
\newtheorem{lemma}{Lemma}
\newtheorem{cor}{Corollary}
\newtheorem{remark}{Remark}
\newtheorem{ex}{Example}
\newtheorem{defi}{Definition}

\usepackage{booktabs}  % nicer lines for tables: \toprule \midrule\bottomrule
\usepackage{mathtools}  % coloneq
\usepackage[normalem]{ulem}
\usepackage[dvipsnames]{xcolor}
\usepackage{placeins}  % FloatBarrier
\usepackage[hidelinks]{hyperref}
\usepackage{upquote}
\usepackage{xtab}

\DeclareMathOperator{\diag}{diag}

\newcommand{\norm}[1]{\left|\left|#1\right|\right|}

\newcommand{\RR}{{\mathbb R}}

\newcommand{\set}[1]{\left\lbrace #1 \right\rbrace}
\let\set\set
\newcommand{\comment}[1]{\relax}
\newcommand{\tbsm}[1]{\left[\!\begin{smallmatrix}#1\end{smallmatrix}\!\right]}
\newtheorem{algorithm}{Algorithm}

\newcommand{\gitlab}[1]         {git\-lab.com/\-#1}

\date{}

\medskip 
\author{
Thomas Mejstrik
\thanks{University of Vienna;  {e-mail: \tt\small
thomas.mejstrik@gmx.at}}\ ,\quad  
Vladimir Yu. Protasov 
\thanks{University of L'Aquila, Italy; {e-mail: \tt\small
vladimir.protasov@univaq.it}}, 
}

\title{Constructive solution of \\ the common invariant cone problem}

\begin{document}
\maketitle

\begin{abstract}

Sets of~$d\times d$ matrices sharing a common invariant cone 
enjoy special properties, which are widely used in applications. 
 However, finding this cone or even proving its existence/non-existence
 is hard.  This problem is known to be  algorithmically undecidable for general sets of matrices. We show that it can nevertheless be efficiently solved in practice. 
 An algorithm that for a given finite set of matrices, either finds a common invariant cone or proves its  non-existence is presented. Numerical results demonstrate that 
it works   for a vast majority of matrix sets. 
 The structure and properties 
 of the minimal and maximal invariant cones are analyzed. Applications to dynamical systems and combinatorics are considered.

\bigskip

\smallskip

\begin{flushright}

\end{flushright}

\end{abstract}
\bigskip

\begin{center}
\large{\textbf{1. Introduction}}
\end{center}
\bigskip

We consider a compact family~$\cA$ of real $d\times d$ matrices and a 
convex closed cone~$K\subset \re^d$ with an apex at the origin. The cone is assumed to be {\em proper}, i.e.,  it is nondegenerate (possesses a nonempty interior) and pointed, i.e., $K \cap (-K) = 0$.   
The cone~$K$ is {\em invariant} for~$\cA$ if for all $A\in \cA$, we have $AK\subset K$. 
The aim of this paper is to obtain a method to  answer the question: 
Does a given family of matrices possess an invariant cone and, if the answer is affirmative, to construct it. The method consists of two algorithms. The first one 
(Primal-Dual Algorithm) proves within finite time the non-existence  of invariant cone. 
If this algorithm applied to a given set~$\cA$ does not halt, then 
we apply the Polyhedral Cone Algorithm which constructs an invariant cone.  

We also show   that if the set of invariant cones of~$\cA$ 
is nonempty, then it possesses the minimal and maximal by inclusion elements. 
An algorithm of constructing the minimal and maximal invariant cones
in an explicit form is derived. Then we address several applications and consider examples. 
\smallskip 

\textbf{1.1. Motivation}. Nonnegative matrices enjoy many special properties which are 
partially described by the Perron-Frobenius theory. Some of them are true for the 
more general case of matrices with an invariant cone. 
 This concerns, for example, 
the existence of Perron eigenvectors. 
According to the Krein-Rutman theorem~\cite{KR}, 
every matrix which 
leaves a cone invariant possesses a Perron eigenvector in this cone, i.e., 
an eigenvector whose eigenvalue is nonnegative and is equal to the spectral  radius of the matrix (all the definitions are given below).

If several matrices share a common invariant cone~$K$, then they act 
``in the same direction'' meaning that they respect the partial order defined by~$K$. 
For every~$\ba$ from the interior of the dual cone, the 
linear functional~$f_{\ba}(\bx) \, = \ (\ba, \bx), \ \bx \in K$, 
is asymptotically equivalent  to the Euclidean norm~$\|\bx\|$. 
Hence, inside the  cone~$K$,  the Euclidean norm 
is equivalent to the value of a certain linear functional and the same is true 
for the operator norms of matrices that leave~$K$ invariant. 
This gives particular spectral and asymptotic properties~\cite{BZ, J12} 
which are put to good use in many  in applications. 
Sets of matrices with common invariant cones appeared  in 
combinatorics and number theory~\cite{GP13, JPB, P00, WZ}, 
neural networks~\cite{EMT},  chaos, dynamical system, and control theory~\cite{BV, DHH, FLS, GS, H}. 
They are applied  in the study of random matrix products~\cite{H, P, P13}
and antinorms~\cite{GZ15, P22}. 

The existence of a common invariant cone simplifies 
the computation of the joint spectral characteristics of matrices  
such as the joint and lower spectral radii~\cite{GP13, JPB, Mej} and the Lyapunov 
exponents~\cite{P, PJ13}.

 \medskip 

 \textbf{1.2. The invariant cone problem. What is known?} A necessary condition for the existence of a common invariant cone 
 is expressed in terms of leading eigenvectors.  For a given matrix~$A$, we denote by~$\rho(A)$ 
 its spectral radius, i.e.,  the maximal modulus  of its eigenvalues.
 \begin{defi}\label{d.5}
  An eigenvalue of a matrix~$A$ is called {\em leading} or~{\em Perron} and denoted by~$\lambda_{\max}(A)$ if it is equal to~$\rho(A)$.
  The corresponding nonzero eigenvector is called leading. 
If~$A$ possesses an invariant cone~$K$, then its leading eigenvectors that belong 
  to~$K$ are called  {\em Perron eigenvectors}.
 \end{defi}
 According to Definition~\ref{d.5}, the notion of Perron eigenvector (in contrast to leading eigenvectors and  
 Perron eigenvalues) depends not only on the matrix~$A$ but also on its 
 invariant cone~$K$.
Only those  leading eigenvectors of~$A$ that belong to~$K$ are Perron.  

By the Krein-Rutman theorem, if all matrices from a set~$\cA$ share a common invariant cone,  
then 
each of them has a Perron eigenvalue. Moreover, all products of those matrices 
and all their linear combinations with nonnegative coefficients have Perron eigenvalues. 
This is a {\em necessary condition} for the existence of a common invariant cone. 
This condition, however, is not sufficient even for one matrix~\cite{V}. 
% A criterion of existence of invariant cone for a given matrix
% was presented in~\cite{V}. 

 \begin{defi}\label{d.7}
 Let~$\cA$ be a compact set of real $d\times d$ matrices. 
 The {\em invariant cone problem} for~$\cA$ is to decide whether all matrices of~$\cA$
 share a common invariant cone. The {\em constructive invariant cone problem} 
 consists in deciding the existence of a common invariant cone and, if it exists, finding at least one  invariant cone. 
 \end{defi}
 The invariant cone problem is usually considered  for finite families of $d\times d$
 matrices~$\cA = \{A_1, \ldots , A_m\}$. The case of one matrix~($m=1$) is completely analysed in~\cite{V}. 
Finite families of  $2\times 2$ matrices are studied in~\cite{EMT, RSS}, 
families of commuting matrices  are considered in~\cite{FV},
families of diagonal matrices are in~\cite{RSS}. 

For general $d\times d$ matrices, a criterion of existence of a common invariant cones 
 was presented in~\cite{P10c}.
It works under the mild assumption that the family~$\cA$ is {\em irreducible},
i.e., the matrices~$A_1, \ldots , A_m$ do not have common nontrivial invariant subspaces. 
The criterion  involves the co-called $L_1$-spectral radius:
$$
\rho_1(\cA) \quad = \quad \lim_{k\to \infty}\, 
\Bigl[ m^{-k}\sum_{i_1, \ldots , i_k}\norm{A_{i_k}\ldots A_{i_1}}\Bigr]^{1/k}.
$$
 See~\cite{P97, W} for more about 
this value. Theorem~1 from~\cite{P10c} asserts that the matrices~$A_1, \ldots , \allowbreak A_m$
share a common invariant cone if and only if~$\rho_1(\cA)$ is equal to 
the largest positive eigenvalue of {\em the mean matrix}~$\, \bar A\, =\, \frac1m \, \sum_{i=1}^m A_i$.
Thus, to decide the existence of a common invariant cone it suffices merely to compare two numbers: $\rho_1(\cA)$ and~$\lambda_{\max}(\bar A)$.  
However, this 
simple criterion is hardly applicable in practice since the computation 
of $L_1$-spectral radius~$\rho_1(\cA)$ can be difficult. Besides, this criterion ensures only the existence of 
an invariant cone and says nothing about its construction. The construction of 
 invariant multicones  was addressed in~\cite{BZ}.

 \medskip 

\textbf{1.3. The complexity of the problem.}  The following result shows that 
the  invariant cone problem is hard. 
\medskip 

\noindent \textbf{Theorem A}~\cite[Theorem~2]{P10c}. 
{\em 1) The invariant cone problem is
algorithmically undecidable for families of four matrices with integral entries.

2) This problem  is NP-hard, for families that consist of four matrices, whose entries take
at most four values: either zero or $\pm a$ or $b$,
where $a$ and $b$ are some integers.}
\medskip 

%Note that Theorem~A concerns only the existence of the invariant cone.  
%The construction of invariant cone  is actually a separate problem. 
This discouraging fact, nevertheless, leaves some hope for efficient solutions
of the invariant cone problem. 
First of all, Theorem~A is proved for matrices 
with a common invariant subspace. We do not know if it holds under the irreducibility assumption, 
 although we believe that it does. 
Second, there may exist algorithms performing well for a large part of matrix families or for some classes of matrices. 
The algorithm derived in Section~5 demonstrates  its efficiency 
for a vast majority of  families of matrices. 
For example, in dimension~$d=5$, by more than~$300$ tests with random 
matrix sets, in~$99\%$ cases  the algorithm solves the problem completely, i.e., 
either finds  a common invariant cone or proves its non-existence. 

 \medskip 

 \textbf{1.4. The direct approach does not work.} 
 Among possible methods  for the invariant cone problem, there is one that seems to be most straightforward and natural. Let us need to verify if a given irreducible  family~$\cA = \{A_1, \ldots , A_m\}$
 possesses a common invariant cone. We construct an increasing nested sequence of cones
 $\{K_{i}\}_{i \in \n}$, where $K_1$ is the Perron 
 eigenvector of the meant matrix~$\bar A$,    
 and for each~$n$, the cone $K_{n+1}$ is the conical hull of $K_n$ and of its images~$A_iK_n, \ i=1, \ldots , m$. 
 If at some iteration we have~$K_{n+1} = K_n$, then~$K_n$ is a common  invariant cone for~$\cA$.  
 This is the idea of the {\em Direct Algorithm} (Algorithm~\ref{alg_direct}). 
It is easily realizable, see Section~4 
for details. 
A quite unexpected observation is that in the numerical experiments  the Direct 
Algorithm never halts  and never produces an invariant cone. Respectively, it cannot decide the 
existence of an invariant cone within finite time.

 This is true even for matrices 
that obviously share an invariant cone, for example, for positive matrices,  
which leave invariant the positive orthant~$\re^d_+$.   
Theorem~\ref{th.30} reveals  theoretical explanations of this phenomenon. 
It asserts that for non-singular matrices, the Direct  Algorithm does not  halt. Nevertheless, this algorithm can be used 
for approximating the minimal common invariant cone, provided that it exists. Moreover, 
it is applied  
as a part of the Primal-Dual Algorithm for deciding the existence of invariant cones (Section~4).

 \smallskip 

\textbf{1.5. The roadmap of the main results}. We begin with proving several properties of 
Perron eigenvectors. Theorem~\ref{th.10} asserts that 
if an irreducible matrix set~$\{A_1, \ldots , A_m\}$ possesses an invariant cone~$K$, then 
the mean matrix~$\bar A$  
has a simple Perron eigenvalue and, respectively, 
a unique,  up to normalization, Perron  
eigenvector~$\bv$. Moreover, in this case either~$\bv$ or~$-\bv$  belongs to~${int}\, K$. 

In Section~3  we prove 
that every compact irreducible
matrix family~$\cA$ whose matrices share an invariant cone 
 possesses minimal and maximal invariant cones. The minimal cone~$K_{\min}$ 
has the following characteristic property: for every invariant cone~$K$ of~$\cA$, 
either~$K$ or~$-K$ contains $K_{\min}$. Similarly, either~$K$ or~$-K$ is 
contained in $K_{\max}$. The existence of~$K_{\min}$ and~$K_{\max}$ is not obvious. 
Moreover, it may not be true for reducible families. Thus, the irreducibility assumption is significant, which can be shown by 
simple counterexamples. 
 We establish special duality relations between~$K_{\min}$ and~$K_{\max}$ 
 and also some basic properties of those cones.
 For example, by the Krein-Rutman theorem,~$K_{\min}$ contains all simple leading eigenvectors of matrices from $\cA$, 
 of all  products of those matrices  and of sums of those products.
 This means that for every simple leading eigenvector~$\bv$
of the aforementioned matrices, either~$\bv$ or~$-\bv$ belongs to~$K_{\min}$.   
 However~$K_{\min}$ does not  coincide with the convex hull of Perron eigenvectors of those  matrices. 
  
 In Section~4 we prove that the Direct Algorithm  converges 
 to~$K_{\min}$, however,  never terminates within finite time unless some of the matrices 
 of~$\cA$ 
 are singular.
Then we present  the Primal-Dual  Algorithm
(Algorithm~\ref{alg_dual})
for deciding the existence of a common invariant cone. This algorithm 
halts if and only if the family does not possess an invariant cone,
(Theorem~\ref{th.35}). Thus,  if the family~$\cA$ does not have an invariant cone, 
then the Primal-Dual  Algorithm proves it  within finite time.
In practice it ensures the  non-existence  within a few iterations.

In Section~5 we derive the Polyhedral Cone Algorithm (Algorithm~\ref{alg_approx}) for constructing an invariant cone
(not necessarily minimal) of a finite family of matrices. In most of numerical experiments 
this algorithm has found an invariant cone, whenever the family  possessed it.  
 Thus Algorithms~\ref{alg_dual} and~\ref{alg_approx}
being applied together provide a compete solution to the constructive invariant cone problem
 for most of matrix families.

Let us stress that the  invariant cone obtained by the Polyhedral Cone Algorithm is never minimal. 
Two problems arise in this context:
1) Find a criterion that ensures that~$K_{\min}$ is polyhedral;
2) If the answer to question 1) is affirmative, construct~$K_{\min}$. 
Theorem~\ref{th.50} solves Problem~1 and finds a structure of the polyhedral cone~$K_{\min}$.
Then Algorithm~\ref{alg_minimal} (the ``Minimal Cone Algorithm'') presented in Section~6 constructs~$K_{\min}$.
If the minimal cone is polyhedral, then  that  algorithm halts and finds it explicitly within finite time (Theorem~\ref{th.60}).
Applying this algorithm to the family of transposed matrices, we obtain the maximal invariant cone~$K_{\max}$. 

The performance  of the presented methods is demonstrated in numerical results in Section~7. 
In Section~8 we consider applications to dynamical systems and combinatorics.

\smallskip

 \textbf{1.6. Notation}. 
Throughout the paper we denote vectors  by bold letters 
and their components by the standard letters.  We use the standard notation~$(\vardot , \vardot)$
for the scalar product in~$\re^d$,  ${ int}\, X$ and~$ \partial\, X$ for  the interior and the boundary, respectively, of a set~$X \subset \re^d$. 

 A convex 
body is a convex compact set with a nonempty interior. A convex cone~$K\subset \re^d$
is a set such that for all~$\bx, \by \in K$ and~$t\ge 0$, we have $\bx+\by \in K, \, t\, \bx \in K$.  
If the converse is not stated, all cones are assumed to be closed, nondegenerate  and pointed.
 A ray starting at the origin is called 
{\em extreme generator} or an~{\em edge} of~$K$ if it is an intersection 
of~$K$ with some hyperplane. 
A ray~$\ell \subset K$ starting at the origin is an edge if and only if it does not contain a midpoint of a segment 
with ends in~$K\setminus \ell$. Any nonzero  vector~$\bv\in \ell$ is also called extremal.  Let~$S$ be a cross-section of~$K$ with a hyperplane~$L$ 
not passing through the origin. Then the extreme points of~$S$ are precisely 
the points of intersection of~$H$ with extreme generators. 
It follows from the Krein-Milman theorem that the cone is a conical hull of its 
extreme generators. A cone  with finitely many edges is polyhedral.

The dual cone to a cone~$K$ 
is~$K^*\, = \, \bigl\{\bx \in \re^d\, : \, \inf_{\by \in K}(\bx, \by) \, \ge \, 0\, \bigr\}$.
The convex hull of a set~$X \subset \re^d$ is denoted by~${co}\, (X)$. 
The conic hull of a set~$V \subset \re^d$ is~${cone} (V)  = {co} \set{\alpha V:\alpha\geq 0}$.

In what follows,~$\cA$ is a compact set of real $d\times d$  matrices, 
$A^T$ denotes the transposed matrix to a matrix~$A$, $\, \cA^T \, = \, 
\{A^T\, : \, A\in \cA\}$. 
By $\cP_{\cA}$ we denote the closure of the set of positive linear  combinations of all 
products of matrices from~$\cA$, where the products are without ordering and with repetitions permitted.
We use the short notation~$\cP = \cP_{\cA}$ and~$\cP^T = \cP_{\cA^T}$.
The diagonal matrix with entries $x_1, \ldots , x_n$ on the diagonal is denoted by $\diag(x_1,\ldots,x_n)$.

\vspace{1cm}

\begin{center}
\large{\textbf{2. Simple Perron eigenvectors}}
\end{center}
\bigskip 

We begin with establishing  one special property of irreducible matrix sets 
with a common invariant cone. It will be used in the next sections.

A Perron eigenvector of a matrix~$A$ is called {\em simple} if 
it is a unique, up to multiplication by a constant,  eigenvector corresponding to the Perron eigenvalue. 
 If a compact family~$\cA$ possesses a common invariant cone~$K$, then every matrix from the set~$\cP_{\cA}$ has a Perron eigenvector in~$K$. 
We are going to show that at least for one matrix from~$\cP_{\cA}$,
 this eigenvector is simple and belongs to the interior of~$K$.  
We begin with the case of finitely many matrices. The arithmetic mean of matrices 
of a finite family~$\cA = \{A_1, \ldots , A_m\}$ will be referred to as the {\em mean matrix}
and denoted as~$\bar A \, = \, \frac{1}{m} \sum_{i=1}^m A_i$. Clearly, if~$\cA$ has an invariant cone, then~$\bar A$  has the same invariant cone.  
\begin{theorem}\label{th.10}
If an irreducible finite family~$\cA$
possesses a common invariant cone~$K$, then  the Perron eigenvector 
of the mean matrix~$\bar A$ is simple and 
belongs to $\, {int}\, K$. 
\end{theorem}
{\tt Proof}. We first show that every Perron eigenvector of~$\bar A$
belongs to the interior of~$K$. Assume the converse:  there  is a  
nonzero vector~$\bv \in \partial K$ such that~$\bar A \bv = \lambda \bv$ for some~$\lambda \ge 0$. Denote by~$\Gamma$ 
the minimal by inclusion face of~$K$ containing~$\bv$. 
There exists a unit vector~$\bn \in \re^d$ such that~$(\bn, \bx) \le 0$
for all~$\bx \in K$ and the equality  holds precisely when~$\bx\in \Gamma$. 
Let us show that~$A_i\Gamma \, \subset \, \Gamma$ for all~$A_i \in \cA$. 
If this is not the case, then there is a point~$\by \in \Gamma$ and~$A_j \in \cA$ such 
that~$A_j\by \notin \Gamma$, i.e., $\, \bigl(\bn, A_j\by\bigr) < 0$. 
Therefore,~$
\bigl(\bn, \bar A\by\bigr) \,
=
\,  \frac1m \, \bigl(\bn, \sum_{i=1}^m A_i\by\bigr) \,
=
\, \frac1m \, \sum_{i=1}^m\, \bigl(\bn,  A_i\by\bigr)\, < \, 0
$.
Since~$\bv$ is an interior point of~$\Gamma$, we have~$\bv - \tau \, \by \in \Gamma$ for all sufficiently small~$\tau >0$. 
On the other hand, 
$\, 0\, =  \, \bigl(\bn, \bar A(\bv - \tau \, \by)\bigr) \, = \, 
\bigl(\bn, \bar A\bv\bigr) \, - \, \tau \bigl(\bn, \bar A\by\bigr) \, = \, 
\bigl(\bn, \lambda \bv\bigr) \, - \, \tau \bigl(\bn, \bar A\by\bigr)\, = \, 
\lambda \bigl(\bn, \bv\bigr) \, - \, \tau \bigl(\bn, \bar A\by\bigr)\, = \, 
 \, - \, \tau \bigl(\bn, \bar A\by\bigr) \, > \, 0$, which is a contradiction. 

Thus, $\bv \, \in \, {int}\, K$. It remains to show that 
this is a unique leading eigenvector of~$\bar A$ up to normalization. 
If, to the contrary, there is 
a  vector~$\bv'\in \re^d$ not collinear to~$\bv$ and such that~$\bar A \bv' = \lambda \bv'$,
 then 
the restriction of the operator~$\bar A$ to the two-dimensional 
subspace~$V \, = \, {span}\, \{\bv, \bv'\}$ is the identity operator multiplied by~$\lambda$. 
Since~$\bv \in {int}\, K$,   it follows that the set~$K\cap V$ is a 
proper cone in~$V$. Every nonzero vector~$\bu$ on the boundary of this cone 
belongs to the boundary of~$K$. On the other hand,~$\bu$ is a Perron eigenvector of~$\bar A$
since~$\bu \in V$. Therefore, $\bu$ must be an interior point of~$K$. 
The  contradiction completes the proof. 

{\hfill $\Box$}
\medskip

\begin{cor}\label{c.10}
Let an irreducible compact matrix family~$\cA$ possess an invariant cone~$K$. 
Then there is a matrix in~$\cP_{\cA}$ whose Perron eigenvector is simple
and belongs to ${int}\, K$. 
\end{cor}
{\tt Proof}. It is well-known that every  irreducible compact set of matrices contains an irreducible 
finite subset. By Theorem~\ref{th.10}, the arithmetic mean of matrices from this subset has a simple Perron eigenvector in~$\, {int}\, K$. 

{\hfill $\Box$}
\medskip 

 Actually the mean matrix in Theorem~\ref{th.10} can be replaced   by any 
positive linear combination of matrices: 
\begin{cor}\label{c.13}
Let an irreducible family~$\cA = \{A_1, \ldots , A_m\}$ possess an invariant cone~$K$
and let~$\{h_i\}_{i=1}^m$ be a set of strictly positive numbers. 
Then the matrix~$\, \sum_{i=1}^m h_iA_i$ has a simple Perron eigenvector that  belongs to ${int}\, K$. 
\end{cor}
{\tt Proof}. Apply Theorem~\ref{th.10} to the family of 
matrices~$\tilde A_i \, = \, \frac{m h_i}{\sum_{i=1}^m  h_i}\, A_i, \quad i=1, \ldots , m$. 

{\hfill $\Box$}
\medskip

\begin{defi}\label{d.8}
Let a matrix~$A$ possess a simple leading eigenvector. Then the leading eigenvectors of~$A$
and~$A^T$ are called, respectively, the right and left Perron eigenvectors of~$A$. 
\end{defi}
By Theorem~\ref{th.10} if an irreducible family~$\cA$ has an invariant 
cone~$K$, then its mean matrix~$\bar A$ has the right and left Perron eigenvectors 
that belong to~${int}\, K$ and~${int}\, K^*$, respectively.

\bigskip 

\begin{center}
\large{\textbf{3. The minimal and maximal invariant cones}}
\end{center}
\bigskip

A matrix family may possess several invariant cones. 
 In this case there should exist a minimal cone being 
the intersection of all invariant cones. This conclusion seems to be obvious, 
 however,  it is actually wrong. Indeed, if~$K$ is an invariant cone, then so is the opposite cone~$-K$. 
On the other hand,~$K \cap (-K) = \{0\}$. This implies that 
the intersection of all invariant cones of a given matrix family is always zero.    
Therefore, the concept of the minimal invariant cone has to be understood 
differently.   
\begin{defi}\label{d.10}
An invariant cone~$K$ of a matrix family~$\cA$ is 
called {\em minimal} if for every invariant cone~$K'$, 
we have either~$K\subset K'$ or~$K\subset (-K')$. 
The {\em maximal} invariant cone is defined similarly: 
for every invariant cone~$K'$, 
we have either~$K'\subset K$ or~$(-K') \subset K$. 
\end{defi}
Let us remember that all the cones are assumed to be proper. 
In what follows we denote the minimal and maximal invariant cones by~$K_{\min}(\cA)$ 
and~$K_{\max}(\cA)$ respectively. Certainly, they are defined only if~$\cA$
has at least one invariant cone. This condition, however, is not sufficient for 
 the existence of~$K_{\min}$ and~$K_{\max}$. 
  For example, the identity matrix has infinitely many invariant cones,  neither of them is minimal or maximal. 
Nevertheless, if the family is irreducible, then 
among its invariant cones there always exist  minimal and maximal ones.  
\begin{theorem}\label{th.20}
If an irreducible matrix family possesses an invariant cone, then 
it has  both the minimal and the maximal invariant cones. Moreover, $\, K_{\max}(\cA) \, = \, 
 [K_{\min}(\cA^T)]^*$.  
\end{theorem}
Clearly, the minimal and the maximal cones are unique. 
The proof of  Theorem~\ref{th.20} is based on the following key lemma: 
\begin{lemma}\label{l.10}
For arbitrary invariant cones~$K_1, K_2$ of an irreducible  matrix family, 
one of the sets~$K_1\cap K_2$ or~$K_1\cap (-K_2)$  is a proper cone 
and the other is zero.    
\end{lemma}
{\tt Proof.} It suffices to consider finite matrix families and then 
generalize as in the proof of Corollary~\ref{c.10}. 
By Theorem~\ref{th.10} applied to the family~$\cA = \{A_1, \ldots , A_m\}$ 
and to its invariant cone~$K_1$, the mean matrix~$\bar A = \frac1m \sum_{i=1}^m A_i$
 possesses a simple Perron eingenvector~$\bv_1 \in {int}\, K_1$. 
Similarly, $\bar A$ has a  simple Perron eingenvector~$\bv_2 \in {int}\, K_2$. 
Therefore, those vectors are collinear and hence, ${int}\, K_1$ intersects the 
interior of either~$K_2$ or of~$-K_2$. Assume it intersects~$ {int}\, K_2$. 
Then~$K_{+} = K_1 \cap K_2$ is a proper invariant cone of~$\cA$. 
If the intersection~$K_{-} = K_1 \cap (-K_2)$ is nonzero, then $K_{-}$ is also an 
invariant cone for~$\cA$. The irreducibility of~$\cA$ implies that $K_{-}$ is a proper cone. 
Thus, we obtain two invariant cones~$K_{+}, K_{-}$ of $\cA$ such that 
$K_{+}\cap K_{-} \, = \, \{0\}$, since those cones are contained in~$K_2$ and~$-K_2$,  respectively. 
However, as it is shown  above, one of those intersections must be a proper cone. 
The contradiction completes the proof.

{\hfill $\Box$}
\medskip

\noindent {\tt Proof of Theorem~\ref{th.20}.} 
Choose one invariant cone~$K$ of the family~$\cA$.  
By Lemma~\ref{l.10}, every other invariant cone~$K'$, 
possibly, after replacing it by its opposite $-K'$, 
intersects~$K$ by a proper cone.  Thus, $K$ intersects all other 
invariant cones of~$\cA$
by proper cones, which are also invariant for~$\cA$. The set of those cones 
is closed with respect to finite intersections. 
Now  the compactness argument proves that 
all those cones  possess a nonzero intersection,  which  is also an invariant cone for~$\cA$. 
It is 
proper due to the irreducibility of~$\cA$. This proves the existence of the minimal cone. 

To obtain the maximal cone we consider the family of transposed matrices~$\cA^T$. 
It possesses the minimal cone~$K_{\min}(\cA^T)$. 
Every other invariant cone of~$\cA^T$ contains either~$K_{\min}(\cA)$ or its opposite.  
Therefore, every invariant cone of~$\cA$ is contained either in~$[K_{\min}(\cA^T)]^*$
or in~$ - \, [K_{\min}(\cA^T)]^*$. 
Hence, the former is the maximal invariant cone of~$\cA$. 

{\hfill $\Box$}
\medskip 

\smallskip

The minimal invariant cone possesses the following self-similarity property: 
it coincides with the conic hull of its images under the action of operators from~$\cA$:  
\begin{prop}\label{p.15}
For every irreducible matrix family, we have~${K_{\min} \, = \, 
{co}\, \{ AK_{\min}\, : \, A\in \cA\}}$. 
%and~$K_{\max} \, = \, \bigcap\limits_{A \in \cA} \, AK_{\max}$.
\end{prop}
{\tt Proof.} The cone~${{ co}\, \{ AK_{\min}\, : \, A\in \cA\}}$ 
is invariant for~$\cA$ and is contained 
in~$K_{\min}$.  Hence, it coincides with~$K_{\min}$. % The proof for~$K_{\max}$ is the same. 

{\hfill $\Box$}
\medskip 

Applying Theorem~\ref{th.10} we immediately obtain 
\begin{cor}\label{c.15}
If an  irreducible matrix family~$\cA =\{A_1, \ldots , A_m\}$
has an invariant cone, then the simple Perron eigenvector of the mean matrix~$\bar A$
 belongs to~${int}\, K_{\min}$. 
\end{cor}

\bigskip 

\begin{center}
\large{\textbf{4. The Direct Algorithm and the Primal-Dual Algorithm}}
\end{center}
\bigskip 

In this section we consider the Direct Algorithm
of constructing invariant cones. We are going to see that 
this  
algorithm in general does not  terminate within finite time and, respectively, 
does not find an invariant cone (and even does not prove its existence),  provided 
that some mild assumptions are satisfied.  Nevertheless, it   
can be modified to the Primal-Dual Algorithm that efficiently 
proves the non-existence of invariant cones. 
\smallskip 

The algorithm begins with choosing 
a simple leading eigenvector~$\bv$ of some matrix~$A \in \cP_{\cA}$. 
By Corollary~\ref{c.15} 
we can always take the  mean matrix~$A = \bar A$. If this matrix does not have a simple leading 
eigenvector, then~$\cA$ does not have an invariant cone.

\begin{algorithm}[Direct Algorithm]\label{alg_direct}
~\begin{flushleft}
{\bfseries Given: } An irreducible family of $d\times d$ matrices 
$\cA = \{A_1, \ldots , A_m\}$.\\
{\bfseries Result upon Termination: } Invariant cone for~$\cA$, or a proof that no cone exists
\rule{0.7\textwidth}{0.8pt}\\
Take an arbitrary matrix from~$\cP_{\cA}$ that has a simple leading eigenvector~$\bv$\\
Set $\cV_0 \coloneq \{\bv\},\ K_0 \, = \, {cone}\, \cV_0$\\
{\bfseries for} $j = 1,2,\ldots$\\
$\qquad$ Set $\cV_{j} \coloneqq \emptyset$\\
$\qquad$ {\bfseries for} $\bx\in {\cA}\cV_{j-1}$ {\bfseries do}\\
%$\qquad\qquad$ Construct some $\bT$ satisfying Definition~\ref{def:tree}\hfill $(\ast)$ \qquad ~\\
$\qquad\qquad$ {\bfseries if}  $\bx \not\in K_{j-1}$\\
$\qquad\qquad\qquad$ Set $\cV_{j} \coloneqq \cV_{j} \cup \bx$ \\
$\qquad$ Set $K_j =  {cone}\, \bigcup\limits_{i=0}^j\cV_i$\\
$\qquad$ {\bfseries if} $K_j = \RR^d$\\
$\qquad$ $\qquad$ {\bfseries return} \emph{``No invariant cone exists''}\\
$\qquad$ {\bfseries elseif} $\cV_{j} = \emptyset$\\
$\qquad$ $\qquad$ {\bfseries return} $K_j$
\end{flushleft}
\end{algorithm}

\begin{remark}\label{r.15}
{\em The Direct Algorithm produces an increasing nested sequence of  cones 
$\{K_{j}\}_{j\ge 0}$. 
These cones become full-dimensional in some $j$th iteration,~$j\le d-1$. 
Indeed, if for some~$i \ge 0$, the linear spans of~$K_{i}$ and~$K_{i+1}$ coincide,
then this span is a common invariant subspace for~$\cA$ and hence, 
this is the entire~$\re^d$. Thus, 
the dimension of~$K_{i}$ strictly increases with~$i$ until it becomes equal to~$d$. 
Hence, $K_{d-1}$ has dimension~$d$. 
}
\end{remark}

We are going to prove that the cones~$K_j$ produced by the Direct Algorithm tend to~$K_{\min}$
as~$j\to \infty$ although do not coincide with it. The convergence means that the sets 
$K_j\cap B $, where $B$ is the unit ball, converge to~$K_{\min}\cap B$ in the 
 Hausdorff metric. The convergence is equivalent to that the closure 
 of the union $\bigcup_{j\ge 0} K_j$ coincides with~$K_{\min}$. This is 
 easily shown by convexity. 

\begin{theorem}\label{th.30}
If a finite irreducible family~$\cA$ has an invariant cone, then, for every initial vector~$\bv$
(a simple leading eigenvector of some matrix from~$\cP_{\cA}$), 
we have~$K_{j} \, \to \, K_{\min}$ as~$j\to \infty$. 

If all matrices from~$\cA$ are non-singular and 
$\bv$ is not an extreme vector of~$K_{\min}$, then the sequence~$\{K_{j}\}_{j\ge 0}$
strictly increases and none of the cones~$K_{j}$ is invariant. In particular, 
none of them coincides with~$K_{\min}$. 
\end{theorem}

\noindent {\tt Proof.} Since
the cone~$K_{\min}$ contains all simple leading eigenvectors of matrices from~$\cP_{\cA}$, 
we have~$\bv \in K_{\min}$ and therefore,~$K_0 \subset K_{\min}$. 
Then by induction it follows that~$K_j \subset K_{\min}$ for all~$j$. 
Furthermore,  the closure of the union of all~$K_{j}$
is an invariant cone of~$\cA$ which is contained in~$K_{\min}$. 
 Hence, it coincides with~$K_{\min}$. 
Thus, $K_{j}\, \to \, K_{\min}$ as~$j\to \infty$. 
If all matrices from~$\cA$ are non-singular and~$\bx$ is a non-extreme vector 
of~$K_{\min}$, then all the vectors~$A\bx, \ A\in \cA$, are also non-extreme. 
This way we show by induction that all the sets~$\cV_j$ produced by the Direct Algorithm 
do not contain  extreme vectors of~$K_{\min}$. 
On the other hand, 
$K_{\min}$ possesses at least one extreme generator. Hence, none of the cones~$K_j$
coincides with~$K_{\min}$ and consequently, none of them is invariant.

{\hfill $\Box$}
\medskip

Invoking Corollary~\ref{c.15} we see that the leading eigenvector~$\bv$ of the mean 
matrix~$\bar A$ lies in  the interior of~$K_{\min}$ and hence, 
is not an extreme vector of this cone. 
\begin{cor}\label{c.25}
If all matrices of the family~$\cA$ are non-singular and 
$\bv$ is the leading eigenvector of the mean matrix~$\bar A$, 
then the Direct Algorithm does not halt. 
\end{cor}

\begin{figure}
\centering
\includegraphics{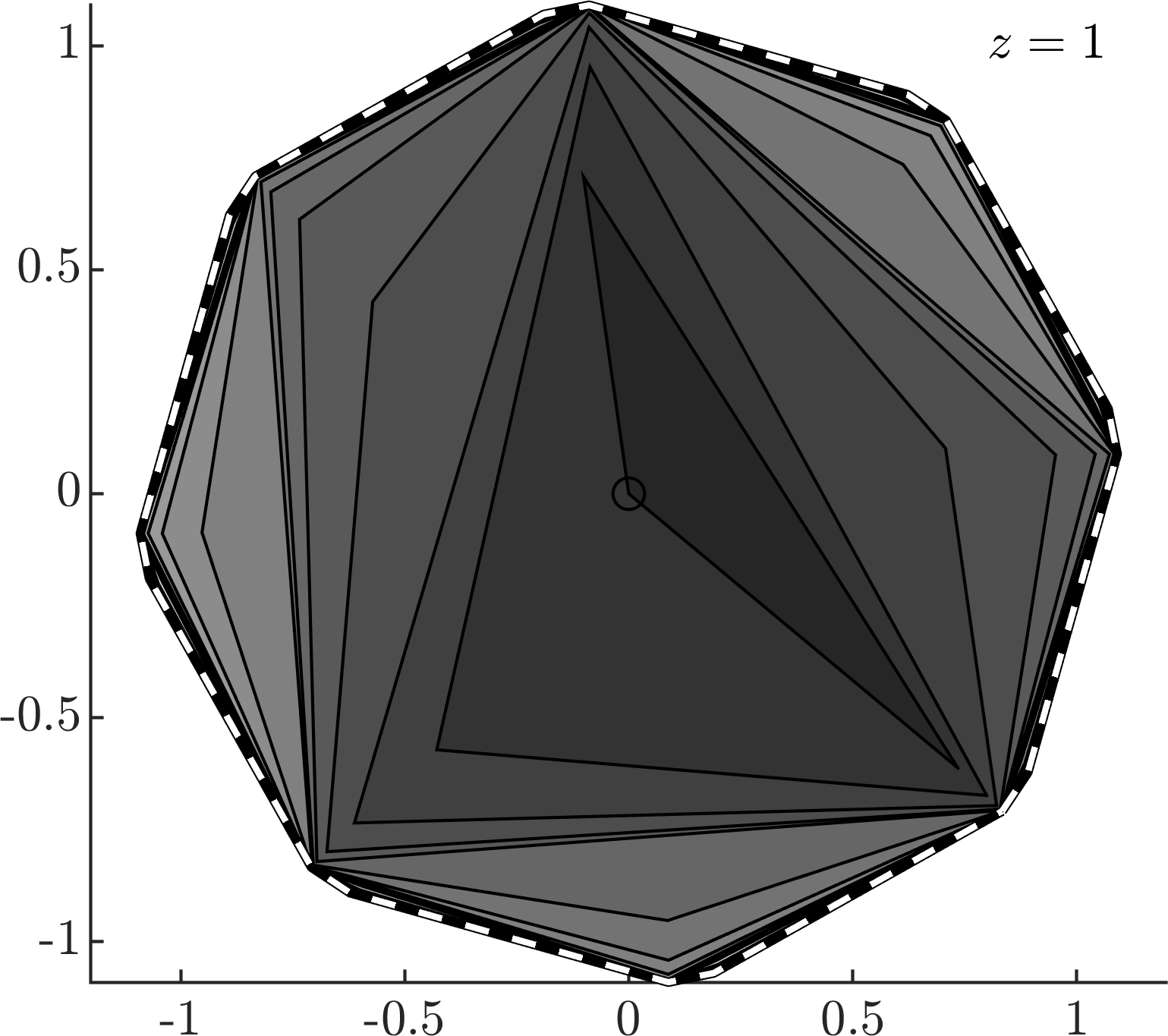}
\caption{{\footnotesize The sequence of cones~$K_j$ generated by the Direct Algorithm for matrices from Example~\ref{ex_direct_algorithm}. The cone~$K_{\min}$ is plotted with a dashed line. The sequence~$K_j$ converges to~$K_{\min}$ but does not reach  it.
}}

\label{fig_direct_algorithm}

% %MATLAB CODE TO GENERATE THE PICTURE
% %DO NOT REMOVE THIS
%clf
%hold on
%K = {center};
%for i = 2:20
%    K{i} = [B{1}*K{i-1} B{2}*K{i-1} K{i-1}];
%    K{i} = normalizematrix( K{i} );
%    nrm = polytopenorm( K{i}, K{i}, 'kone', 'output','ub', 'v',0 );
%    K{i}(:,nrm<1) = [];
%end
%
%for i = numel( K ):-1:3
%    plotm( K{i}, 'Kone', 'center',center, 'hull', 'patch',1, 'color',[i i i]/numel(K) )
%end
%axis equal
%%alpha .5
%
%plotm( center, 'Kone', 'center',center, 'ok', 'ms',5 )
%plotm( [nfo_B.blockvar{1}.cyclictree.VV{:}], 'Kone', 'hull', 'center',center, 'k-', 'linewidth',1.5 )
%plotm( [nfo_B.blockvar{1}.cyclictree.VV{:}], 'Kone', 'hull', 'center',center, 'w--', 'linewidth',1 )
%%plotm( [nfo_B.blockvar{1}.cyclictree.v0{:}], 'Kone', 'center',center, 'k.', 'ms',10 )
%%plotm( [nfo_B.blockvar{1}.cyclictree.v0{:}], 'Kone', 'center',center, 'ko', 'ms',5 )
%figure( 1 )
%
%exportfig.AXISE( [-1.2 1.2] )
%exportfig.TICKS()
%exportfig.LABEL()
%exportfig.SZE(2,2)
%exportfig.TEXT( '$z=1$', .8, 1 )
%exportfig.SAVE( 'direct_algorithm', [], '-png' )
\end{figure}

\begin{ex}\label{ex_direct_algorithm}
{\em Let 
\begin{equation*}
\cA = \set{%
\begin{bmatrix}2&3&6\\4&1&8\\0&0&14\end{bmatrix},\ %
\begin{bmatrix}-1&-1&0\\\phantom{-}1&-1&0\\\phantom{-}0&\phantom{-}0&\sqrt{2}\end{bmatrix}
}
\end{equation*}
Fig.~\ref{fig_direct_algorithm} shows the minimal invariant cone~$K_{\min}$ and the sequence of cones~$K_j$ generated by the Direct Algorithm. Every polygone is the intersection of the corresponding cone with the plane~$L = \{(x,y,z) \in \re^3: \ z=1\}$. 
 The starting one-dimensional cone~$K_0 = \{\bv\}$ is marked with a circle in the middle of the picture. The cone~$K_1$ is shown by the black triangle, the cone~$K_2$ is slightly lighter, 
 $K_3$ is still lighter, etc. The cone~$K_{\min}$ is the octagon plotted with a dashed line. We see that 
the sequence of cones~$K_j$ converges to~$K_{\min}$ but does not reach  it. 
The picture generated using \texttt{plotm} from \texttt{ttoolboxes}~\cite{ttoolboxes}).}
\end{ex}

\begin{remark}\label{r.35}. 
{\em Note that for different initial vectors, the Direct Algorithm produces 
different sequences of cones~$\{K_j\}_{j\ge 0}$.  By Theorem~\ref{th.30}
all those sequences converge to~$K_{\min}$, provided that~$\cA$ has an invariant cone. 
 If~$\cA$ consists of non-singular  matrices, then for 
 all initial vectors, except for, possibly, extreme vectors  of~$K_{\min}$, 
 those sequences do not contain any invariant cone of~$\cA$. 

}
\end{remark}

In spite the negative conclusion of Theorem~\ref{th.30}, the Direct Algorithm can nevertheless be put to good use
as an auxiliary procedure to any method of the invariant cone problem. Being applied 
to the transpose family~$\cA^T$, it gives a stopping criterion. This procedure will be referred to as {\em the dual complement}.  
\medskip

\textbf{The dual complement}. Let us have an iterative algorithm of constructing an invariant cone,
which produces an embedded sequence of cones~$\{K_j\}_{j \ge 0}$. If some of them is invariant, then the algorithm halts. 
It will be referred to as the {\em Primal Algorithm}. Note that this an arbitrary algorithm, not necessarily the Direct Algorithm.
  \medskip 

Take the smallest~$q$ for which the cone~${ int} K_q\, \ne \, \emptyset$.  
Take the leading eigenvector~$\bv^*$ of the matrix~$\bar A^T \, = \, 
\frac1m \sum_{i=1}^m A_i^T $ and normalize it so that~$(\bv^*, K_q) \ge 0$.  We use the notation~$(X,Y) \, = \, \{(\bx,\by)\, : \ \bx\in X, \, \by\in Y\}$.

If this  is impossible, then STOP: 
$\cA$ does not have an invariant cone.

Start the Direct Algorithm for the family~$\cA^T$ producing the 
sets~$\cV_j^*, \, j\ge 0$. Denote~$\cU_j^* \, = \, \cup_{s=0}^{j}\cV_s$. In each iteration we 
check that~$(K_{q+j}\, , \, \cU_{j}^*)\, \ge \, 0$. 
If this is not the case, then STOP: 
$\cA$ does not have an invariant cone. 
Indeed, since by Theorem~\ref{th.20}~$K_{\max} =  [K_{\min}(\cA^T)]^*$ (possibly, 
after replacing~$K_{\max}$ by~$- K_{\max}$) 
and~$\cU_{j}^* \subset K_{\min}(\cA^T)$, we see that 
$(K_{\max}, \cU_{j}^*) \, \ge \, 0$ and hence,~$(K_j, \cU_{j}^*) \, \ge \, 0$. 

 \bigskip
 
The dual complement can be applied to every Primal Algorithm which 
constructs an embedded sequence 
of cones. In particular, choosing the Primal Algorithm 
 to be the Direct Algorithm 
starting with the initial leading 
eigenvector~$\bv$ of the mean matrix~$\bar A$, we obtain the following method:

\smallskip

\begin{algorithm}[Primal-Dual Algorithm]\label{alg_dual}

{\em Initialization}.  We have a finite irreducible matrix family~$\cA$ of~$d\times d$ matrices.  
Take, respectively,  the  right and left Perron 
eigenvectors~$\bv, \, \bv^*$ of the mean matrix~$\bar A$ such 
that~$(\bv, \bv^*) \geq 0$. If at least one of them does not exist or is not simple, 
then STOP: $\cA$ has no invariant cones.  

Define the sets~$\, \cU_0 \, = \, \cV_0 \coloneq \set{\bv}$, 
$\, \cU_0^* \, = \, \cV_0^* \coloneq  \set{\bv^*}, \ K_0 = {cone}\, \cU_0, \ 
K_0^* = {cone}\, \cU_0^*$.
%Set $V_{max}$ (maximum number of vertices) to some value.

\smallskip

{\em $j^{th}$ Iteration step}. Set~$\cV_j = \emptyset$. For each 
 $\bx\in {\cA}\cV_{j-1}$ we set~$\cV_{j} \coloneqq \cV_{j} \cup \bx$ if~$\bx \notin K_{j-1}$, 
 otherwise, we leave~$\cV_{j}$ as it is. When all~$\bx$ are exhausted, we 
 set~$\cU_{j}  \coloneqq \cU_{j-1} \cup  \cV_{j}$. 
\smallskip 
Similarly, set~$\cV_j^* = \emptyset$ and  for each 
 $\bx^*\in {\cA}^*\cV_{j-1}^*$ we 
 set~$\cV_{j}^* \coloneqq \cV_{j}^* \cup \bx$ if~$\bx^* \notin K_{j-1}^*$, 
 otherwise, we leave~$\cV_{j}^*$ as it is. When all~$\bx^*$ are exhausted, we 
 set~$\cU_{j}^*  \coloneqq \cU_{j-1}^* \cup  \cV_{j}^*$. 
\smallskip

Compute the signs $\Sigma_j$ of $(\bx^\ast, \bx)$ among all 
pairs $\bx^\ast\in\cU^\ast_{j}$, $\bx\in\cV_{j}$,  and all pairs 
$\bx^\ast\in\cV^\ast_{j}$, $\bx\in\cU_{j}$.

\smallskip 

{\em Stopping criterion}. If there is at least one negative sign in $\Sigma_j$, then $\cA$ does not have an invariant cone. 
%If $\abs{\cV} > V_{max}$ then the algorithm could not decide whether no invariant cone exists.
\end{algorithm}

In most of numerical examples the Primal-Dual Algorithm proves the non-existence of invariant cones within a few iterations. The following theorem asserts that if $\cA$ does not have in invariant cone,
then the Primal-Dual Algorithm always terminates within finite time.

\begin{theorem}\label{th.35}
The  Primal-Dual  Algorithm  applied to an irreducible finite matrix set~$\cA$
halts if and only if~$\cA$ does not have an invariant cone. 
\end{theorem}

This theorem is a straightforward corollary  of the following result: 

\begin{theorem}\label{th.40}
Each of the following conditions is equivalent to the non-existence of 
invariant cone for an irreducible family~$\cA$: 
\smallskip 

\noindent \textbf{a)} For some~$j\in \n$, 
 the set~$\cU_j\setminus \{0\}$ contains~$d+1$ vertices 
of a simplex that covers the origin;
\smallskip 

\noindent  \textbf{b)} For some~$j\in \n$, the sets~$\cU_j$ and~$ \cU_j^*$
 contain two nonzero vectors (one in each set) forming an obtuse angle.  
 \smallskip 
 
\noindent  \textbf{c)} The Primal-Dual Algorithm applied to~$\cA$ halts. 
\end{theorem}
\noindent {\tt Proof.} Denote~$\cU \, = \, \bigcup_{j\ge 0}\cU_j$ and
~$ \cU^* \, = \, \bigcup_{j\ge 0} \cU_j^*$. 
Obviously, b) and c) are equivalent. Let us show that 
a) is equivalent to the non-existence of an invariant cone. Then we do the same with b). 
% It is also obvious that a) implies b) and implies the non-existence of an invariant cone. 
% It remains to prove that b) implies a) and that b) implies the non-existence of an invariant cone. 

Obviously, a) implies the non-existence, since the set $\cU_j$ cannot be contained in 
a proper cone.  Let us establish the converse. 
To every point $\ba \in \re^d\setminus \{0\}$ we 
associate its polar~$[\ba]^*\, = \, \bigl\{\bx \in \re^d: \ (\bx, \ba) \le 0\bigr\}$.
Assume none of the simplices with vertices in~$\cU_j\setminus \{0\}$
contain zero.  Then every~$d+1$ half-spaces from the 
set~$\bigl\{[\ba]^*: \ \ba \in \cU \setminus\{0\}\,  \bigr\}$ have a nonempty intersection.
By Helly's theorem~\cite{MT}, this implies that all of them have a nonempty 
intersection. Denote it by~$Q$.
This is a cone  invariant with respect to all matrices from~$\cA^T$. 
The irreducibility implies that~$Q$ is non-degenerate. 
 Hence~$Q^*$ is a pointed cone invariant with respect to  all matrices from~$\cA$. 
 Applying again the irreducibility argument we conclude that~$Q^*$ is nondenerate. 
 Thus, $\cA$ possesses a proper invariant cone, which is a contradiction.  
 \smallskip  

It now remains to show that b) is equivalent to the non-existence of an invariant cone. 
We begin with necessity. Since the set~$\cU_{\, d-1}^*$  
has a full-dimension, so is~$\,  \cU^*$. 
Therefore, if $\, (\cU, \cU^*) \ge 0$, then 
the conical hull of~$\cU$, which is an invariant cone for~$\cA$,
is pointed. Thus, $\cA$ has a proper invariant cone.  
  \smallskip  
 
Conversely, assume~$\cA$ possesses an invariant cone.  Then 
  every point~$\bx \in \cU_j$ belongs to~$K_{\min}(\cA)$ and hence, to~$K_{\max}(\cA)$.  
  Respectively, every point~$\bx^* \in  \cU_j^*$
  belongs to~$K_{\min}(\cA^T) =   [K_{\max}(\cA)]^*$. Therefore, 
  ${(\bx^*\, , \, \bx) \ge 0}$.

{\hfill $\Box$}
\medskip

\begin{center}
\large{\textbf{5. Construction of a common invariant cone}}
\end{center}
\bigskip 

Any algorithm that solves the constructive invariant cone problem for a given matrix family~$\cA$
has to realize two options: 1) Prove the non-existence of common invariant cones
in the case when~$\cA$ does not have them; 2)  construct at least one invariant cone otherwise.  
\smallskip 
  
 The part 1 is done by the 
Primal-Dual  Algorithm (Algorithm~\ref{alg_dual}). 
 In this section we realize part 2. By Theorem~\ref{th.30}, this cannot be done by the Direct Algorithm because the  polyhedral cones of the produced embedded family~$\{K_j\}_{j\in \n}$ are ``too narrow'': they are all contained in~$\, {int}\, K_{\min}$. Hence, to obtain an invariant cone, one needs to enlarge those cones~$K_j$ in each iteration. This is the idea of the 
 Polyhedral Cone Algorithm. The enlargement will be realized by 
 the ``scaling'' of the cones with respect to a ray generated by a fixed vector~$\bc \in \, {int}\, K_{\min}$.  We begin with auxiliary results.  
 \bigskip 
 
 \begin{center}
\textbf{5.1. The mean matrices whose left and right Perron eigenvectors coincide}
\end{center}
\smallskip 
 
 \begin{lemma}\label{l.20}
 Let~$\cA$ be a set of matrices, $\bv^*, \, \bv$ be, respectively, the 
 left and right Perron eigenvectors of the mean matrix~$\bar A$, 
 $(\bv, \bv^*) > 0$. 
 
 Suppose~$V$ is a nonsingular matrix such that~$\ V\, V^T \bv^*\, = \, \bv$. 
 Then the  left and right eigenvectors of the matrix~$V^{-1}\bar AV$ coincide and are equal to~$\bc =  V^{-1}\bv$. 
 \end{lemma}
 {\tt Proof}. Multiplying both parts of the inequality 
~$\ V\, V^T \bv^* \, = \, \bv$  from the left by~$V^{-1}$, we obtain
 $\ V^{T}\bv^* \, = \, V^{-1} \bv \, = \, \bc$ . 
 Let~$\lambda$ be the Perron eigenvalue of~$\bar A$. Then
$$
\left\{
\begin{array}{lclclcl}
V^{-1} \bar A V\bc &= & V^{-1} \bar A \bv & = & \lambda V^{-1} \bv & = & \lambda \bc \\ 
V^{T} \bar A^T V^{-T}\bc &= & V^{T} \bar A^T \bv^* & = & \lambda V^{T} \bv^* & = & \lambda \bc 
\end{array}
\right. 
$$ 
 Thus, 
 $$ 
\left\{
\begin{array}{lcl}
V^{-1} A V\bc &= &  \lambda \bc \\ 
\bigl( V^{-1} A V \bigr)^T \bc &= & \ \lambda \bc 
\end{array}
\right. 
$$

{\hfill $\Box$}
\medskip

\begin{remark}\label{r.60}
{\em The matrix~$V$ in~Lemma~\ref{l.20} can be computed as a square root 
of a positive definite matrix~$P$ such that~$\bv = P \bv^*$. The latter  can be found
by an explicit formula with an arbitrary parameter~$\kappa > 0$:  
\begin{equation}\label{eq.P}
P \  = \ \frac{\bv\cdot \bv^T}{(\bv^*, \bv)} \ + \ 
\kappa\, \left( I \ - \frac{\bv^*\cdot {\bv^*}^T}{(\bv^*, \bv^*)} \right).
\end{equation}
Then it remains to find the matrix~$V$ such that~$V\, V^T\, = \, P$.  

}
\end{remark}

Passing to the new basis with the transfer matrix~$V$ we obtain the matrix 
family~$\tilde \cA\, =\, V^{-1}\cA V$ with the mean matrix~$V^{-1}\bar A V$. 
The left and right Perron eigenvectors of this matrix coincide and equal to~$\bc = 
V^{-1}\bv$. If~$\cA'$ has an invariant cone~$K'$, 
then~$\cA$ has an invariant cone~$K = VK'$.  
\bigskip 

 \medskip 
 
 \begin{center}
\textbf{5.2. The cone-norm and the matrix scaling}
\end{center}
\bigskip

\begin{defi}\label{def_Et}
Let $K\subseteq\RR^d$ be a proper cone and $\, \bc\, \in \, {int} \, (K\cap K^*)$, 
$\|\bc\|=1$.
Let~$C$ be an orthonormal matrix with the first column $\bc$.
For arbitrary $t>0$,  we define the matrix
\begin{equation}
E_{t,\bc} \ = \  E_t \ =
\ C \, \diag(1,\ t, \ldots, t) \, C^T \ =  \ 
C \, 
\begin{bmatrix}
    1 &   &        & 0 &  \\
      & t &        &   &  \\
      &   & \ddots &   &  \\
    0 &   &        & t &
\end{bmatrix}
\, C^T \ \in \ \RR^{d \times d}\, . 
\end{equation}
\end{defi}

\begin{remark}\label{r.30}
{\em The condition~$\bc \in {int}\, K^*$ is equivalent to~$\bc^{\perp} \cap K \, = \, 
\{0\}$ and to the following property: the cross-section of~$K$ by an arbitrary affine hyperplane orthogonal to~$\bc$ is a bounded set. 
Indeed, if for some~$\ba \in {int}\, K$, the cross-section of the cone~$K$ by the 
hyperplane $\, \ba \, + \, \bc^{\perp}$
is not bounded, then it contains a ray~$\{\ba + t \bh: \ t\in \re_+\}$, where~$\bh \in \bc^T$. 
Therefore, $t^{-1}\ba + \bh \, \in \, K$ for all~$t\ge 0$. 
Therefore, $\bh \in K$, which is a contradiction, since~$\bc^{\perp} \cap K \, = \, 
\{0\}.$ Thus, 
the condition~$\bc \, \in \, {int}\, (K\cap K^*)$ is equivalent to 
that~$\bc$ is an interior point of~$K$ and the cross-section passing through~$\bc$
and parallel to~$\bc^{\perp}$ is bounded. }
\end{remark}

The matrix $E_t$ can be used to ``scale'' cones, as the following lemma demonstrates.

\begin{lemma}\label{thm_Et}
%With $E_t$ as the multiplication, $K$ becomes nearly a vector space (except for the missing inverse).
\begin{enumerate}
\item \label{Et_indepence_of_nullspace} The matrix $E_t$ does not depend on the choice how the matrix $C$ is complemented to a full basis. 
\item \label{Et_multiplicativity} $E_{s} E_{t} = E_{st}$. 
%\item \label{Et_identity} $E_1 = I$ is the identity matrix
%\item \label{Et_0} $E_0=cc^T$ is the orthogonal projection onto $c$ (and thus $E_0 w \neq 0\in\RR^d$ in general)
%\item \label{Et_inf_vector} $\lim_{t\rightarrow\infty} E_t v \in c^{\perp}$ for all $v\in\RR^d$

\item \label{Et_inclusion} $E_{s,\bc} K \subseteq E_{t,\bc} K \Leftrightarrow s \leq t$
for every cone $K$ for which $\bc$ fulfills the properties in Definition~\ref{def_Et}.
%\item \label{Et_inf_K} 
%$
%\lim_{t\rightarrow\infty} E_t K =
%\set{\alpha_v  v + \sum_{w\in c^\perp} \alpha_w w : \alpha_v, \alpha_w \geq 0}
%$,
%i.e. the half-space defined by $c$

\end{enumerate}
\end{lemma}
{\tt Proof.}
Item $1$ follows from the fact that the product 
  $\bc\, \bc^T$
is the orthogonal projection onto~$\bc$.  
Item $2$ follows since $C$ is orthonormal. 
Item $3$ follows from the convexity of $K$.

{\hfill $\Box$}

\medskip 

 The scaling has a simple geometric interpretation. 
Denote by~$\ell$ the ray~$\{\lambda \bc: \ \lambda \ge 0\}$. The operator~$E_t$ respects all affine hyperplans orthogonal to~$\ell$, on  each of them~$E_t$ defines 
the homothety with coefficient~$t$ with respect to the point of intersection of this 
hyperplane with~$\ell$.

%\begin{proof}
%\begin{enumerate}
%\item[\ref{Et_indepence_of_nullspace}.]
%\begin{align*}
%\tbmatrix{c&c^\perp}\tbmatrix{1 & 0\\0 & tI}\tbmatrix{c^T\\{c^\perp}^T} =
%cc^T + t{c^\perp}{c^\perp}^T
%\end{align*}
%and the product ${c^\perp}{c^\perp}^T$ is the orthogonal projection onto $c^\perp$.
%\item[\ref{Et_multiplicativity}.]
%\begin{align*}
%E_{s,c} E_{t,c}
%& = V \operatorname{diag}(1,\ s, \ldots, s) V^T V \operatorname{diag}(1,\ t, \ldots, t) V^T \\
%& = V \operatorname{diag}(1,\ s, \ldots, s) \operatorname{diag}(1,\ t, \ldots, t) V^T\\
%& = V \operatorname{diag}(1,\ st, \ldots, st) V^T \\
%& = E_{st,c}
%\end{align*}
%\item[\ref{Et_identity}.] $E_1 = I$ since $C$ is orthonormal.
%\item[\ref{Et_0}.]
%$E_0 = cc^T$, and since
%% and since $E_0^2 = cc^Tcc^T = cc^T = E_0$ it follows that $E_0$ is a projection.
%$(cc^Tx,y)=(c^Tx,c^Ty)=(x,cc^Ty)$ it follows that $E_0$ is an orthogonal projection.
%$c$ is an eigenvector of $E_0$ corresponding to eigenvalue 1, and thus the projection is onto $\left<v\right>$.
%\end{enumerate}
%\end{proof}

\begin{figure}
\centering
\includegraphics{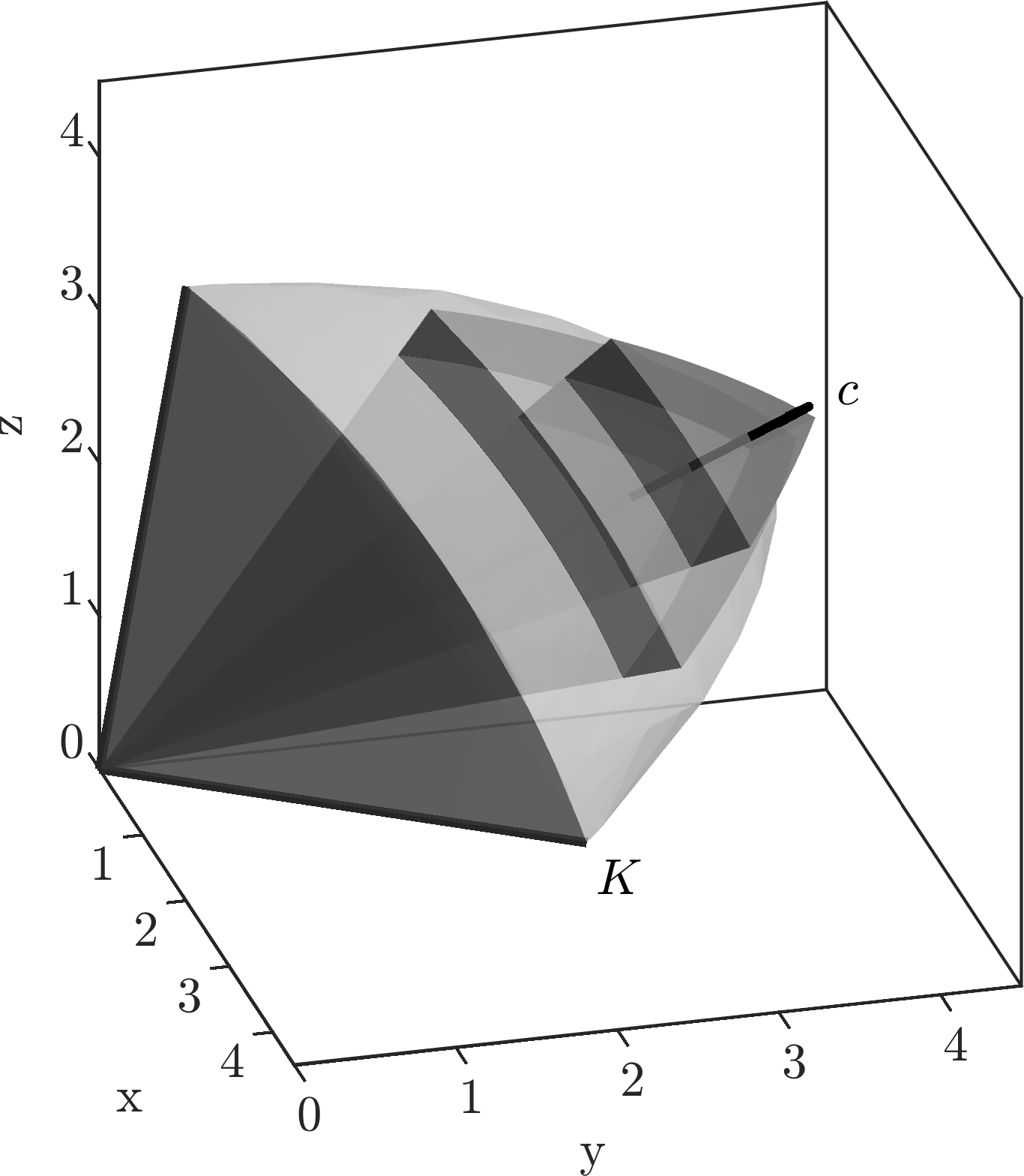}
\caption[A cone generated by three vectors
(the most light one),
  together with its scaled versions by
  factors of $1/2$ and $1/4$ (the darker ones).
  The center vector used for scaling
  is printed as a black line
  ]{%
{\footnotesize   A cone generated by vectors
$\set{%
\tbsm{2\\0\\4},
\tbsm{1\\4\\2},
\tbsm{4\\2\\1}
}$ (the most light one),
  together with its scaled versions by
  factors of $1/2$ and $1/4$ (the darker ones) with respect to the .
   center vector $\bc$.
  (Picture generated using \texttt{plotm} from
  \texttt{ttoolboxes}~\cite{ttoolboxes}).}
  }

\label{fig_Et}
%%%
%% MATLAB CODE TO GENERATE PICTURE
%% DO NOT REMOVE THIS
%c = [4 4 4]';
%c = c/norm(c);
%K = [2 0 4;1 4 2;4 2 1]';
%mu = 2.5;
%clf; hold on;
%
%plotm( c*1.6, 'k', 'arrow',mu*1.5, 'linewidth',1.5 )
%plotm( reshape( [zeros(size(K)); normalizematrix(K)*mu*1.8; nan(size(K))], 3, [] ),  'color',0.15*ones(1,3), 'linewidth',1.5, '-' )
%
%t = .25;
%plotm( Et_kone( t, c ) * K, 'funct',{'kone',mu*2.2},'hull', 'color',0.3*ones(1,3) )
%
%t = .5;
%plotm( Et_kone( t, c ) * K, 'funct',{'kone',mu*2}, 'hull', 'color',0.4*ones(1,3) )
%
%t = 1;
%plotm( Et_kone( t, c ) * K, 'funct',{'kone',mu*1.8}, 'hull', 'color',0.5*ones(1,3) )
%
%axis equal
%
%plotm( 'scene','default' )
%alpha .7
%view( 75, 24 )
%exportfig.AXIS( [0 4.5 0 4.5 0 4.5] )
%exportfig.TICKS(1:1:4, 0:1:4, 0:1:4)
%exportfig.LABEL('x','y','z')
%exportfig.TEXT( '$K$', 4, 2, 0.75 )
%%exportfig.TEXT( '$E_{1/2}K$', 4.25, 2.5, 1.75 )
%%exportfig.TEXT( '$E_{1/4}K$', 4.25, 3, 2.75 )
%exportfig.TEXT( '$c$', 3.6, 3.6, 3.6 )
%box on
%
%exportfig.SZE(2,2)
%exportfig.POST
%%%
%%exportfig.SAVE( 'Et', [], '-png' )
\end{figure}

\begin{defi}
The cone-norm $\norm{\vardot}_K = \norm{\vardot}_{K,\bc}$ is defined on a given cone~$K$
by a vector~$\, \bc\, \in \, {int} \, (K\cap K^*), \, 
\|\bc\|=1$, as follows: 
\begin{equation}
\norm{\bx}_K \ = \ \inf \, \bigl\{ t>0 \, : \ \bx \in E_{t,\bc} K \, \bigr\}
\end{equation}
\end{defi}
On every cross-section~$G$ of the cone~$K$ by a hyperplane orthogonal to 
the ray~$\ell = \{\lambda \bc : \ \lambda \ge 0\}$, the cone-norm is equal to the 
Minkowski functional of the set~$G$ with respect to the point~$G\cap \ell$. 
Since~$G$ is a bounded set~(Remark~\ref{r.30}), it follows  that {\em the conic 
norm is well-defined at every point}. 
\begin{lemma}\label{thm_kone_norm}
The cone-norm $\norm{\vardot}_K \, = \, \norm{\vardot}_{K,\bc}$ fulfills the following properties for all $\bx\in\RR^d$ with $(\bx,\bc)>0$
\begin{enumerate}
\item \label{vecnorm_0} $\norm{\bx}_K \, = \, 0 \Leftrightarrow \bx\ =\ \lambda \bc$ for some $\lambda \geq 0$

\item \label{vecnorm_scalarmult} $\norm{\bx}_K = \norm{\lambda \bx}_K \ = \ 
\norm{\bx}_{\lambda K}$ for all $\lambda > 0$

\item \label{vecnorm_Etmult} $\ t\norm{\bx}_K = \norm{E_t \bx}_K \ = 
\ \norm{\bx}_{E_{t^{-1}} K}$ for all $t > 0$

%\item If $\norm{x}_K \leq \norm{y}_K$, then $\norm{x}_K \leq \norm{x+y}_K \leq \norm{y}_K$ \XX{proof missing}

\item \label{vecnorm_triangle} $\norm{\bx+\by}_K \leq \norm{\bx}_K + \norm{\by}_K$

\item \label{vecnorm_lp}
If $K$ is given as a conic hull of finitely many vertices $\cV\subseteq\RR^d$, then
\begin{equation*}
\norm{\bx}^{-1}_K = \begin{cases}
    \operatorname{maximize} \ t \ \\
    \operatorname{ subject \, to: }\\
    \diag(1,\ t, \ldots, t) \, \tilde{\bx} \ - \ \sum_{\tilde{\bv}\in \tilde{\cV}} \, \alpha_{\tilde{\bv}} \, \tilde{\bv} \ = \ 0\\
    \alpha_{\tilde{\bv}} \, \geq \, 0 \\
    t \ \geq \ 0
\end{cases}
\end{equation*}
where $\tilde{\bx} = C^T \bx$, $\tilde{\cV} = C^T \cV$.
\end{enumerate}
\end{lemma}
{\tt Proof.}
Items~$1$, $2$, and~$3$ follow directly from Lemma~\ref{thm_Et}.
Point~$4$ follows from the convexity of~$K$.
Point~$5$ is straightforward to check.
%\begin{enumerate}
%\item[\ref{vecnorm_0}.]
%$E_0 K = \set{\lambda v:\lambda \geq 0}$.
%If $x=c$ then $\norm{p}_K = \inf \set{ t>0 : x\in E_t K} = 0$.
%If $\norm{p}_K = 0$ it follows that $x\in E_t K$ for all $t>0$, and thus $x=\lambda c$ for some $\lambda \geq 0$.
%
%\item[\ref{vecnorm_scalarmult}.]
%This follows because $K = \lambda K$ for all $\lambda > 0$.
%
%\item[\ref{vecnorm_Etmult}.]
%\begin{align*}
%\norm{x}_K 
% & = \inf \set{ s>0 : x\in E_{s} K} \\
% & = \inf \set{ s>0 : x\in E_{t^{-1}}E_tE_{s} K} \\
% & = \inf \set{ s>0 : E_t x\in E_{st} K} \\
% & = \inf \set{ \tilde{s}/t>0 : E_t x\in E_{\tilde{s}} K} \\
% & = \frac{1}{t} \norm{E_t x}_K
%\end{align*}
%
%The second equality follows similarly.
%%\begin{align*}
%%t\norm{x}_K 
%% & = \inf \set{ s>0 : E_t x\in E_{s} K} \\
%% & = \inf \set{ s>0 : x\in E_{t^{-1}}E_{s} K} \\
%% & = \inf \set{ s>0 : x\in E_{s}E_{t^{-1}} K} \\
%% & = \norm{x}_{E_{t^{-1}}K}
%%\end{align*}
%
%\item[\ref{vecnorm_triangle}.]
%Let $t > \norm{x}_K $, $s > \norm{y}_K$.
%Then, because $K$ is convex,
%\begin{align*}
%E_{\frac{1}{t+s}}(x+y)
%=
%E_{\frac{t}{t+s}}E_{\frac{1}{t}} x
%+
%E_{\frac{s}{t+s}}E_{\frac{1}{s}} y
%\in
%K
%\end{align*}
%Now let 
%$t-\norm{x}_K = \varepsilon$,
%$s-\norm{y}_K = \varepsilon$,
%and it follows that 
%$
%\norm{x+y}
%\leq 
%t + s
%=
%\norm{x}_K + \norm{y}_K
%+ 2\varepsilon
%$.
%Since $\varepsilon$ is arbitrary, the claim follows.
%
%\item[\ref{vecnorm_lp}.]
%This is straight forward.
%
%\end{enumerate}

{\hfill $\Box$}

\newpage 
 
 \begin{center}
\textbf{5.3. Polyhedral Cone Algorithm}
\end{center}
\bigskip

\begin{algorithm}[Polyhedral Cone Algorithm]
\label{alg_approx}
~\begin{flushleft}
{\bfseries Given: } A  set of matrices~$\cA = \{A_1, \ldots , A_m\}$, a number~$t > 1$\\
{\bfseries Result upon Termination: } Invariant cone for~$\cA$
%\Vhrulefill\\
\rule{0.7\textwidth}{0.8pt}\\
{\bfseries Initialization:}\\
%{\bfseries Algorithm:}\Tstrut\\ 
%Search for s.m.p.-candidates $A_1,\ldots, \Pi_M$, set $\lambda \coloneqq \rho(\Pi_1)^{1/\len\Pi_1}$\\
%Scale matrices $\tilde{\cA} \coloneqq \{\lambda^{-1} A_j:j=1,\ldots,J\}$\\
%Compute approximate minimal dual cone $K^\ast$\\
Set $\bar{A} \coloneqq \sum_{i=1}^m  A_i$\\
Compute the leading eigenvectors $\bv$ and $\bv^*$ of $\bar{A}$ and $\bar{A^T}$ respectively\\
Normalize $\bv^*$ so that $(\bv^*, \bv) = 1$\\
Compute positive definite matrix $P$ such that  $\bv \, = \, P\bv^*$ 
(formula~(\ref{eq.P}))\\
Compute a $d\times d$ matrix~$V$ such that $VV^T \, = \, P$\\
Compute the center vector $\bc  \coloneq V^{-1}\bv $\\
Scale matrices: set $\tilde{\cA} \coloneq \{ E_{t,\bc} V^{-1} A_j  V\}$ \\
{\bfseries Main Loop:}\\
Set $\cV_0 \coloneqq \set{\bv_0}$\\
{\bfseries for} $j = 1,2,\ldots$\\
$\qquad$ Set $\cV_{j} \coloneqq \emptyset$\\
$\qquad$ Set $N_{j} \coloneqq 0$\\
$\qquad$ {\bfseries for} $\bx\in \tilde{\cA}\cV_{j-1}$ {\bfseries do}\\
%$\qquad\qquad$ Construct some $\bT$ satisfying Definition~\ref{def:tree}\hfill $(\ast)$ \qquad ~\\
$\qquad\qquad$ {\bfseries if}  $n_{\bx} \coloneqq \norm{\bx}_{\tilde K_{j-1},\bc} > N_j$\\ 
$\qquad\qquad\qquad$ Set $N_{j} \coloneqq \max\{N_{j},\ n_{\bx}\}$\\
$\qquad\qquad\qquad$ Set $\cV_j\coloneqq \cV_j\cup \bx$ \\
$\qquad$ Set $\tilde K_j \coloneqq {cone}\, \bigcup_{n=0}^j \cV_n$\\
$\qquad$ {\bfseries if} $\tilde K_j = \re^d$\\
$\qquad\qquad$ {\bfseries return} \emph{``Invariant cone cannot be computed''}\\
$\qquad$ {\bfseries else if} $N_j \leq t$\\
$\qquad\qquad$ {\bfseries return} cone $K_j = V\tilde K_j$, center $V\bc$
\end{flushleft}
\end{algorithm}

\begin{remark}
{\em If Algorithm~\ref{alg_approx} terminates with \emph{``Invariant cone could not be computed''}
it does not mean that there does not exist an invariant cone for $\cA$.
It merely means that the set of scaled matrices $\tilde{\cA}$ does not possess an invariant cone. The algorithm may terminate for some~$\tilde t < t$. 

The Direct Algorithm (Section~4)  is actually Algorithm~\ref{alg_approx} for~$t=1$. 
In this case, by Corollary~\ref{c.25}, 
if all the matrices are non-singular, then Algorithm~\ref{alg_approx} never terminates.  That is why, we  set in the Initialization~$t> 1$.  Algorithm~\ref{alg_approx}
may work differently for different values of~$t$.}
 
\end{remark}

\begin{theorem}\label{th.70}
If Algorithm~\ref{alg_approx} terminates with a cone~$K_j$,
then~$K_j$ is a common invariant cone for~$\cA$.
\end{theorem}
{\tt Proof.}
By construction, each point $\bv \in \bigcup_{n=0}^j \cV_n$ fulfills the condition
%\begin{equation*}
$\|\tilde{A} \bv\|_{\tilde K_j,\bc} \quad \leq \quad t\, , \ \tilde{A} \in \tilde \cA$.
%\end{equation*}
Therefore, 
\begin{equation*}
\norm{V^{-1}AV \bv}_{\tilde K_j,\bc} \ = \ \frac{1}{t}\, \|\tilde A \bv \|_{\tilde K_j,\bc} \ \leq \ \frac{t}{t} \ = \ 1\, , \qquad A\in \cA\, . 
\end{equation*}%
We see that~$\tilde K_j$ is an invariant cone of the family~$V^{-1}\cA V$, hence, 
the cone $K_j = V\tilde K_j $ is invariant for~$\cA$. 

{\hfill $\Box$}
\medskip

\begin{center}
\large{\textbf{6. When is the minimal cone polyhedral?}}
\end{center}
\bigskip

Algorithm~\ref{alg_approx} constructs a polyhedral invariant cone. 
This cone, however, 
  may not be minimal. Moreover, $K_{\min}$ may not be polyhedral. 
  Two problems arise in this context: 
 \begin{enumerate}
 \item {\em How to find the minimal invariant cone?}
 \item {\em When is the minimal invariant cone polyhedral?}
 \end{enumerate}
  
By Theorem~\ref{th.30}, the Direct Algorithm
gives an approximation of~$K_{\min}$, but does not find it and even does not prove that 
the family has an invariant cone. 
In this section we derive  a method that decides whether the minimal cone is polyhedral and, if 
the answer is affirmative,   finds it.
The method is based on Theorem~\ref{th.50}, which 
establishes  the structure of the minimal cone. 
To formulate it we need to introduce some more notation.  

{\em The cyclic tree of matrix products}. 
Consider a  finite matrix family  $\cA \, = \, \{A_1, \ldots
,  A_m\}$. To an arbitrary  product  $ \Pi \, = \, 
A_{d_n}\cdots  A_{d_1}, \, n\ge 1$ (all products are with repetitions permitted) 
we associate
the {\em cyclic tree} $\cT \, = \, \cT(\Pi)$ defined as
follows. The root is formed by the cycle~$\cR = \cR(\Pi) = (\bv_1, \ldots , \bv_n)$, 
where~$\bv_1$ is the leading eigenvector 
of~$\Pi$ and~$\bv_{j+1}\, = \, A_{d_j}\bv_j, \ j=1, \ldots , n-1$. 
 Thus, $\bv_j$ is the leading eigenvector of 
the product~$A_{d_{j-1}}\cdots A_{d_{1}} A_{d_n} \cdots A_{d_j}$. 
We see that the cycle consists of leading eigenvectors of the cyclic permutations 
of~$\Pi$. The  elements of~$\cR$ are, by definition, the nodes of
zero level.  For every $i\le n$,  an edge (all edges are directed)
goes from $\bv_i$ to $\bv_{i+1}$, where we set $\bv_{n+1} = \bv_1$. At each
node~$\bv_j$ of the root $m-1$ edges start to nodes of the first level.
Those are the nodes~$A_s\bv_j, \, s\ne j$. 
So, there are $n(m-1)$ different nodes on the first level. The
sequel is by induction: from every node~$\bv$ of the~$(k-1)$st level 
$m$ edges go to its ``children'' $A_s\bv, \, s=1, \ldots , m$. 
Thus, the $k$th level contains 
 $n(m-1)m^{k-1}$ nodes. A subset of nodes~$\cS$ of~$\cT\setminus \cR$ is 
 called a {\em minimal cut set} if it has exactly one common node with every  
infinite path along  the tree starting from the root. All the paths are without backtracking. 
A node~$\bv$ is {\em blocked}
by~$\cS$ if the infinite path starting at~$\bv$ intersects~$\cS$. In particular, 
all elements of the cyclic root are blocked. A {\em finite cyclic three} is a subtree of~$\cT$
that consists of nodes blocked by some minimal cut set. Every finite tree contains the cyclic root~$\cR$. We see that every product 
$\Pi$ generates a unique cyclic tree (up to multiplication of~$\bv_1$ by a constant) and infinitely many finite cyclic trees. 

In what follows we identify a nonzero vector~$\bv$ with the ray~${cone} (\bv)\, = \, 
\{\, t\bv \, : \ t\ge \}$. So, the nodes of cyclic trees can be considered as vectors or as rays.

\begin{theorem}\label{th.50}
If the minimal invariant cone  is polyhedral, then the 
set of its edges  is a union of several finite cyclic trees. 
\end{theorem}
Thus, if~$K_{\min}$ is a polyhedral cone, then 
there are several products of matrices from~$\cA$ 
and finite trees generated by those products such that 
the set of all edges of~$K_{\min}$ is precisely the 
elements of those trees. Therefore, the  minimal cone 
can be obtained merely by exhaustion of the matrix products and 
constructing cyclic trees similarly to the Direct Algorithm. This idea is realized in 
Algorithm~\ref{alg_minimal} presented below. 

\smallskip 

\noindent {\tt Proof of Theorem~\ref{th.50}}. 
Choose an arbitrary  edge~$\ell_1$ of~$K_{\min}$. 
It is well known that every extreme point of a convex hull of some sets is an extreme point of one of those sets. Since~$K_{\min}$ coincides with the convex hull of the 
images~$AK_{\min}$ over all~$A\in \cA$ (Proposition~\ref{p.15}), it follows that 
one of the sets~$AK_{\min}$ contains~$\ell_1$.  
On the other hand, $AK_{\min}$ is a polyhedral cone contained in~$K_{\min}$.  Therefore, 
$\ell_1$ is an edge of~$AK_{\min}$. This means that   
$K_{\min}$ has an edge~$\ell_2$ such that~$A\ell_2 = \ell_1$. 
Thus,   {\em every edge of~$K_{\min}$ is an image of another edge 
by means of some operator~$A\in \cA$}. Iterating we obtain a 
chain of edges~$\ell_1\to \ell_2 \to \ldots $ such that~$\ell_j = A(j)\ell_{j+1}$, 
where~$A(j) \in \cA$, $j\ge 1$. Since~$K_{\min}$ has finitely many edges, this chain must cycle: $\ell_{r} = \ell_{r+n}$ for some 
$r$ and~$n \ge 1$. Then~$\ell_r$ corresponds to the Perron eigenvector~$\bv_r$ of the 
product~$A(r)\cdots A(r+n-1)$ of length~$n$, and the cycle of this product~$\{\bv_{r+n-1}, \ldots , \bv_{r}\}$ generates a cyclic 
tree~$\cT_1$. The path along~$\cT_1$ from this cycle to the edge~$\ell_1$ 
consists of edges of~$K_{\min}$. Then we take an arbitrary new edge of~$K_{\min}$ 
(call it~$\ell_{r+n}$), 
which  does not belong  to the set~$\{\ell_1, \ldots , \ell_{r+n-1}\}$. 
 Arguing  as above we build a path~$\ell_{r+n} \to \ell_{r+n+1} \to \ldots $. If it meets
  the previous set of edges, then~$\ell_{r+n} \in \cT_1$, otherwise this path cycles 
  and forms a new root~$\cT_2$, etc.  
  After several iterations the finite set of edges of~$K_{\min}$ will be split into several finite cyclic trees.

{\hfill $\Box$}
\medskip

 Theorem~\ref{th.50} suggests the following method of constructing the minimal cone. 
Make a certain order in the set of products and add them one by one 
constructing the corresponding cyclic trees.  By the same principle as in the Direct Algorithm, 
we remove all newly appeared vertices 
if they are dead (in the conic hull of 
the previously constructed vertices). 
When  in some iteration no new vertices appear, then~$K_{\min}$ is polyhedral 
 and is spanned by the vectors of the constructed cyclic trees. Now we present the formal procedure.   

% new: 

\begin{remark}\label{rem_minmal}
{\em In the algorithms description, we define
 $\min (A,B) \, = \, \min\limits_{a\in A, b\in B} (a,b)$

Algorithm~\ref{alg_minimal} uses the function $\pi_\cA$, where
$\pi[i]$ is the i$^{th}$ in lexical order element,
which is not a power of a shorter one,
in the set $\bigcup_{n\in\n} \{1,\ldots,m\}^n$.
$\pi_\cA[i]$ is the product with factors of $\cA$ given by $\pi[i]$.
E.g.\ for $m=2$ we have: $\pi_\cA[1] = A_1$,
$\pi_\cA[2] = A_2$,
$\pi_\cA[3] = A_2A_1$, $\pi_\cA[4]=A_2A_1^2$,
$\pi_\cA[5]=A_2^2A_1$,
$\pi_\cA[6]=A_2A_1^3$, $\pi_\cA[7]=A_2^2A_1^2$, \dots.
}
\end{remark}

%endnew

\begin{algorithm}[Minimal Cone Algorithm]\label{alg_minimal}
~\begin{flushleft}
{\bfseries Given: } Matrices $\cA$\\
{\bfseries Result upon Termination: } The minimal invariant cone for~$\cA$
%\Vhrulefill\\
\rule{0.7\textwidth}{0.8pt}\\
Run several iterations of the Direct Algorithm for the family~$\cA^T$ and obtain a 
nondegenerate cone~$K^* \subset K_{\min}(\cA^T)$. \\
Set~$\cV_0 \coloneqq  \emptyset$\\ 
{\bfseries for} $j=1,2,\ldots$\\
$\qquad$ Set $\ \cR_j \coloneqq$ the set of leading eigenvectors of $\pi_\cA[j]$\\
$\qquad$ {\bfseries if}  $\ \min\, (K^*\, , \, \cR_j)\, \le \, 0\, $, then\\
$\qquad \qquad$ {\bfseries return} ``No  invariant cone exists''\\
$\qquad$  $\cV_j \coloneqq \set{ \bv\in \cR_j : \bv\notin \, {cone}\, (\bigcup_{n=1}^{j-1} \cV_n) }$\\
$\qquad$ {\bfseries if} $\ \cV_j = \cR_j$\\ 
$\qquad$ {\bfseries for} $\ n=1,\ldots,j$ $\quad$ // $(\ast)$\\
% $\qquad \qquad$ {\bfseries for} $\ \bv\in V_j$\\
% $\qquad \qquad \qquad$ {\bfseries if} \, $\ {cone}\, (\bigcup_{n=1}^{j-1} \cV_n)$\\
% $\qquad \qquad \qquad \qquad$ $\cV_n \coloneq \emptyset$\\
% $\qquad \qquad \qquad \qquad$ $V_n \coloneq \emptyset$\\
$\qquad \qquad$ {\bfseries for} $\bv\in\cA \cV_n$\\
$\qquad \qquad \qquad$ {\bfseries if} $\bv\notin \, {cone}\, (\bigcup_{n=1}^j \cV_n)$\\
$\qquad \qquad \qquad \qquad$ $\cV_n \coloneqq \cV_n \cup \{\bv\}$\\
$\qquad$ {\bfseries if} $\ \min (\,\bigcup_{n=1}^j \cV_n\, ,\, K^\ast\,)\, <\, 0$\\
$\qquad \qquad$ {\bfseries return} ``No  invariant cone exists''\\
$\qquad$ {\bfseries if} no vertices were added in loop $(\ast)$\\
$\qquad\qquad$ {\bfseries return} $K_j \, = \, {cone}\, \bigcup_{n=1}^j \cV_n$
\end{flushleft}
\end{algorithm}

 \begin{theorem}\label{th.60}
 Assume that for 
  an irreducible finite family~$\cA$, the Minimal Cone Algorithm
  halts for some~$j \in \n$. 
  Then $\cA$ has an invariant cone, the minimal cone~$K_{\min}(\cA)$ 
  is polyhedral and is spanned by the 
vectors of the images of the leading eigenvectors from $V_i$.
Conversely, if 
the minimal invariant cone of~$\cA$ is polyhedral, then the algorithm halts, 
provided that all the leading eigenvectors are simple. 
\end{theorem}

{\tt Proof}. {\em Sufficiency.} If the algorithm halts, then  by the construction, 
the cone~$K$ spanned by the obtained cyclic trees is invariant. 
Since all the leading eigenvectors of the products from~$\cP_{\cA}$ are simple, it follows that 
they belong to~$K_{\min}$. Therefore, $K=K_{\min}$ and this cone is polyhedral. 

{\em Necessity.} If~$K_{\min}(\cA)$ is polyhedral, then by Theorem~\ref{th.50}, 
it is spanned by finitely many cyclic trees. 
Let~$\{\cR_i\}_{i=1}^k$ be their  roots and~$r$ be their maximal length. Each~$\cR_i$ consists of extreme vectors of~$K_{\min}$ and 
hence, none of them is in the conical hull of the other roots~$\{\cR_s\}_{s\ne i}$. 
%Therefore, none of~$\cR_i$ will be removed at the step~2 of the algorithm. 
This implies that  the obtained set of products contains all~$\cR_i, \, i=1, \ldots , k$. 
Since the finite trees with the roots~$\cR_i$ generate~$K_{\min}$, it follows that 
after some $j$th iteration we obtain~$K_j = K_{\min}$. 
All other vectors involved in the algorithm belong to~$K_{\min}$, hence, no new 
vectors will appear after the $j$th iteration and the algorithm halts.  
 
{\hfill $\Box$}
\medskip

\begin{remark}\label{r.50}{\em In all steps of the Minimal Cone Algorithm 
we check that some vectors belong to conical hulls of other vectors.
In practice, to account for rounding errors and inaccuracies
arising from floating-point computations, we strengthen this condition and check whether
vectors belong to the interior of the conical hull. 
The algorithm still generates~$K_{\min}$ although some redundant  
cycles and vectors may appear, which do not affect the final result.  
}
\end{remark}
\bigskip 

\begin{center}
\large{\textbf{7. Numerical results}}
\end{center}
\bigskip

\begin{center}
\textbf{7.1 Solution of the invariant cone problem}
\end{center}
\medskip

We present numerical solutions of the constructive invariant cone problem for pairs (${m=2}$) of 
random $d\times d$ matrices with normally distributed entries. We leave 
only those pairs~$\cA = \{A_1, A_2\}$ where both matrices have Perron  eigenvalues, all other pairs are not considered. 
To each pair,  we apply the Primal-Dual Algorithm. If it proves the non-existence of invariant cone, then the problem is solved. Otherwise, we apply 
the Polyhedral Cone algorithm and construct a common invariant cone.
In Table~1 the first column ({\em dim})
is the dimension, the second column is the number of  tests. 
We see that for~$d=2,3$, the problem is completely solved in all 
tests, while for $d=4,5$ there are rare cases ($2\%$ and~$1\%$ respectively) when 
the algorithms could neither construct an invariant  cone nor prove its non-existence.  
  For~$d\ge 6$, the Polyhedral Cone algorithms usually works too long, and we report only the results of the Primal-Dual Algorithm, which works well even in high dimensions.

%\begin{figure}
%\centering
$$
\begin{tabular}{|c|c|c|c|c|c|}
\toprule
% \multicolumn{5}{c}{Gaussian - first orthant} \\\midrule
 dim & \#Tests &  $\cA$ has no invariant cone   & invariant cone of $\cA$ is constructed   &    Unknown     \\\midrule 
   2 &     114 &    \phantom{0}37 (   32\%) &     77 (   68\%) &     0           \\
   3 &     243 &   156 (   64\%) &     87 (   36\%) &     0           \\
   4 &     352 &   276 (   78\%) &     70 (   20\%) &    6 (    2\%) \\
   5 &     331 &   284 (   86\%) &   42 (   13\%) &    5 (    1\%) \\
%   6 &     415 &   329 (   79\%) &     19 (    5\%) &   58 (   14\%) \\
%   7 &     460 &   391 (   85\%) &     10 (    2\%) &   45 (   10\%) \\
   8 &     456 &   421 (   92\%) &      &  35 (    8\%) \\
%   9 &     450 &   412 (   91\%) &      &   34 (    8\%) \\
%  10 &     450 &   424 (   94\%) &                &   24 (    6\%) \\
%  11 &     450 &   413 (   91\%) &           &   33 (    8\%) \\
  12 &     450 &   418 (   92\%) &          &   32 (    8\%) \\
%  13 &     450 &   415 (   92\%) &            &   30 (    7\%) \\
%  14 &     450 &   426 (   94\%) &         &   21 (    5\%) \\
%  15 &     450 &   416 (   92\%) &       &   30 (    7\%) \\
%  16 &     450 &   421 (   93\%) &           &   29 (    7\%) \\
  17 &     450 &   436 (   96\%) &                   &   14 (    4\%) \\
\bottomrule
\end{tabular}
$$
\smallskip 

\noindent {\footnotesize \textbf{Table~1.} Solution of the constructive invariant cone problem for pairs of random $d\times d$ matrices with positive leading eigenvalues}

\medskip 
%The example was computed using the script \emph{./demo/kone/\-demo\_ass\-ess\_kone\_algo\-rithm%\_pri\-mal\_ran\-dom\_matri\-ces\_b.m}
%at Git commit \emph{9d933e7}.

In the next experiment we take a pair of random matrices~$\{A_1, A_2\}$ with uniformly distributed values in 
$[0,\,  1]$ and find the maximal~$\lambda > 0$, for which the 
family~$\cA_\lambda = \{A_i -\lambda I, \, i=1,2\}$ has an invariant cone. 
We determine via bisection
the lower and upper bounds $(\lambda_{-},\ \lambda_{+})$ for~$\lambda$. 

In Table~2 we report for each dimension~$d$,
the the mean length~$\mu$ of the difference 
$\lambda_{+} - \lambda_{-}$
and its standard deviation $\sigma$. 
For the cases where an invariant cone~$K$ could be constructed, we report 
the number $\# V$ of extreme generators  of~$K$, 
again via its mean and standard deviation.
%The example was computed using the script
%\emph{./demo/kone/\-demo\_ass\-ess\_kone\_algo\-rithm\_primal\_dual.m}
%at Git commit \emph{6b96aa1}.

$$
\begin{tabular}{|c|cc|cc|}
    \toprule
   dim   &  $\mu_d$   & $ \sigma_d$ & $\mu_{\# V}$ & $\sigma_{\# V}$ \\ \midrule

   %3    & $0.008176$ & $0.016677$  &   $18228$    &     $16880$     \\ % 29 Tests, very old commit
   %3    & $0.006561$ & $0.007749$  &    $3127$    &     $1225$      \\ %  4 Tests, rho ~= 1
    3    & $0.008002$ & $0.012360$  &    $3199$    &     $4584$      \\ %  9 Tests, rho == 1
   %4    & $0.016479$ & $0.012657$  &    $5064$    &     $4016$      \\ %  5 Tests, rho ~= 1
    4    & $0.044718$ & $0.081986$  &    $5855$    &     $5532$      \\ % 15 Tests, rho == 1
   %5    & $0.143229$ & $0.120750$  &    $3360$    &     $1514$      \\ %  3 Tests, rho ~= 1
    5    & $0.024231$ & $0.020743$  &    $7858$    &     $5541$      \\ % 10 Tests, rho == 1
   %6    & $0.150757$ & $0.131910$  &    $9469$    &     $4921$      \\ %  3 Tests, rho ~= 1
    6    & $0.099854$ & $0.179621$  &   $12506$    &     $5025$      \\ % 10 Tests, rho == 1
    7    & $0.037231$ & $0.011333$  &   $14384$    &     $3075$      \\ % 10 Tests, rho == 1
\bottomrule 
\end{tabular}
$$
{\footnotesize \textbf{Table~2.} Solution of the constructive invariant cone problem
for matrices of the form~$A - \lambda I$.}

\bigskip

\begin{center}
\textbf{7.2 Finding the minimal invariant cone}
\end{center}
\bigskip

%\begin{ex}\label{ex.30}{\em
We applied Algorithm~\ref{alg_minimal} for around~$7000$ pairs of Boolean 
(with entries~$\{0,1\}$) $3\times 3$ matrices. In $87\%$ cases Algorithm~\ref{alg_minimal} constructs the minimal invariant cone.
This means, in particular,  that at least for $87\%$ cases, the minimal cone is polyhedral.
Note that all those pairs obviously have the invariant cone~$\re^3_+$, but their minimal cone is different. 
 
 The set
$
\cA\, =\, \left\{
     \begin{bsmallmatrix}
     1 & 1 & 0 \\ 
     1 & 1 & 1 \\ 
     1 & 0 & 1
     \end{bsmallmatrix}
     ,\, 
     \begin{bsmallmatrix}
     1 & 1 & 1 \\ 
     0 & 0 & 1 \\ 
     1 & 0 & 0
     \end{bsmallmatrix}
\right\}
$
had the longest starting product, namely $A_1^{2}A_2^{2}A_1A_2^{2}A_1A_2^{2}$.
In this case the minimal invariant cone is generated by one 
cyclic tree with the root of length~$10$. 
For the family
$
\cA\, =\, \left\{
     \begin{bsmallmatrix}
     0 & 1 & 1 \\ 
     1 & 0 & 0 \\ 
     1 & 0 & 0
     \end{bsmallmatrix}
     ,\, 
     \begin{bsmallmatrix}
     1 & 1 & 1 \\ 
     1 & 0 & 1 \\ 
     0 & 0 & 0
     \end{bsmallmatrix}
\right\}\, , 
$
 the invariant cone is generated by $7$ cyclic trees (this is the maximal number in all 
 the experiments) corresponding to the products: 
$\{
A_1,\,\allowbreak 
A_1^{3}A_2,\,\allowbreak 
A_1^{2}A_2,\,\allowbreak 
A_1^{2}A_2^{2},\,\allowbreak 
A_1A_2,\,\allowbreak 
A_1A_2^{2},\,\allowbreak 
A_2
\}
$.
%}
%\end{ex}

\begin{remark}\label{r.80}
{\em Our numerical implementations is incorporated into the \verb|ttoolboxes|~\cite{ttoolboxes}.
The experiments are conducted using an AMD Ryzen 3600@3.6 GHz, 6 cores, 64 GB RAM, Matlab R2023b. The results presented in Table~1 were computed using the script
 \emph{./demo/kone/\-demo\_ass\-ess\_kone\_algo\-rithm\_pri\-mal\_ran\-dom\_matri\-ces\_b.m}
at Git commit \emph{9d933e7}. The results of Table~2 were computed 
by the script \emph{./demo/kone/\-demo\_ass\-ess\_kone\_algo\-rithm\_pri\-mal\_dual.m}
at Git commit \emph{6b96aa1}. 

}
\end{remark}

\bigskip 

\begin{center}
\large{\textbf{8. Applications}}
\end{center}
\bigskip

Among applications of the invariant cone problem we spot several ones that are closely related to each other. They concern  the joint spectral characteristics of matrices and  
linear switching systems. 
\medskip

\noindent \textbf{Linear switching systems}. A {\em discrete time linear switching system} is 
the difference equation~$\bx(k+1) \, = \, A(k)  \bx(k), \, k\ge 0; \ \bx(0)\, = \, \bx_0$, 
where~$\bx(k) \in \re^d$ and  for each~$k$, the $d\times d$ matrix~$A(k)$ belongs to a given compact set of matrices~$\cA$. If all matrices of~$\cA$ are nonnegative, then the system is called 
{\em positive.} Positive systems are,  as a rule, much easier to 
analyze~\cite{BCV, FoV, JS}. 
The same holds for a more general class of systems with a common invariant cone, 
see~\cite{GP13}.
Having found an invariant cone for~$\cA$ one can apply methods for positive systems. 
Below we consider two joint spectral characteristics of the set~$\cA$
that are  used in the study of linear switching systems. 
\smallskip 

\noindent \textbf{The lower spectral radius} of a compact set of matrices~$\cA$
is 
$$
\check \rho(\cA)\, = \, \lim_{k\to \infty}\min_{A(i) \in \cA}\|A(k)\ldots A(1)\|^{1/k}\, .
$$
 This concept has many applications in functional analysis, dynamical systems, and 
 combinatorics~\cite{BCV, DK, DST, GP13, G, JPB, P00}. For linear switching systems, 
 the  lower spectral radius is responsible for the {\em stabilisability}.  
 The  system is {\em stabilisable} when 
there is a sequence~$\{A(k)\}_{k\ge 0}$ such that for every~$\bx_0$, we have 
$\bx(k) \to 0$, see~\cite{BCV, FoV, JS}. 
 The stabilisability  is equivalent to that~$\check \rho(\cA) < 1$~\cite{G}. 
It is known that the computation of~$\check \rho(\cA)$ for general 
finite set of matrices is an algorithmically undecidable problem~\cite{TB}. Moreover, 
this quantity 
depends  discontinuously on the  matrices from~$\cA$~\cite{BM}.  However, for matrices 
with a common invariant cone, the situation is different. There is an efficient {\em Invariant Polytope Algorithm} that for many families, finds the precise value of the lower spectral radius~\cite{GP13, GZ15, Mej}. 
 \smallskip

 \smallskip 

\noindent \textbf{The multiplicative Lyapunov exponent}
is the 
limit~$\lambda  =  \lim_{k\to \infty} \frac1k  E\, \log \, 
\|A(k)\ldots A(1)\|$, where~$A(\cdot)$ are independent and identically  distributed 
random matrices, $E$ is the expectation, see~\cite{DP, H, LS} for more details. For simplicity, we assume that~$A(\cdot)$
are uniformly distributed on a given finite set~$\cA$. The computation 
of~$\lambda$ is hard and this quantity may   discontinuously depend 
on the matrices. On the other hand, for nonnegative matrices,  
this dependence is continuous and there are methods for estimating~$\lambda$ with a 
given precision~\cite{Key, P, PJ13}. Most of them can be applied to matrices with a common  invariant cone without any change.    

 \smallskip 

\noindent \textbf{The mortality problem} for a finite family of matrices~$\cA$
is to decide the existence of zero product (with repetitions permitted). 
It is known to be algorithmically undecidable for general integer matrices and NP-hard 
for $0-1$ matrices~\cite{TB, Pet}.  Mortality implies that~${\check \rho(\cA) = 0}$
(for rational matrices those properties are equivalent).  
Hence, for many families with a common invariant cone, the mortality can be efficiently decided. 
\bigskip

We choose one simple  example that illustrates all the aforementioned applications. 
\bigskip 

\begin{ex}\label{ex.10}
{\em Let $\cA = \{ A_1, A_2 \}$ with
\begin{equation}
A_1 = \begin{bmatrix}  \phantom{-}0 & 1 & 1\\ \phantom{-}0 & 1 & 1\\-1 & 1 & 0\end{bmatrix},\quad 
A_2 = \begin{bmatrix}   -1 & 0 & 0\\ \phantom{-}0 & 1 &1 \\ \phantom{-}0 & 1 & 0\end{bmatrix}
\end{equation}
The matrices are not positive, nevertheless, they do have a common invariant cone. 
Applying Algorithm~\ref{alg_minimal} we  compute  
$K_{\min}(\cA) = {cone}\, \{k_1,k_2,k_3,k_4\}$ with
\begin{equation*}
k_1=\begin{bmatrix}1\\1\\1\end{bmatrix},\ %
k_2=\begin{bmatrix}1\\1\\0\end{bmatrix},\ %
k_3=\begin{bmatrix}-1\\\phantom{-}2\\\phantom{-}1\end{bmatrix},\ %
k_4=\begin{bmatrix}-1\\\phantom{-}1\\\phantom{-}1\end{bmatrix},
\end{equation*}
Indeed, simplifying the notation~$K_{\min}(\cA) = K$, we have 
\begin{equation*}
A_1 K =
{cone}\, 
\left\{
\begin{bmatrix}1\\1\\0\end{bmatrix},
\begin{bmatrix}1\\1\\0\end{bmatrix},
\begin{bmatrix}1\\1\\1\end{bmatrix},
\begin{bmatrix}1\\1\\1\end{bmatrix}
\right\}, \quad
A_2 K =
{cone}\, 
\left\{
\begin{bmatrix}-1\\\phantom{-}2\\\phantom{-}1\end{bmatrix},
\begin{bmatrix}-1\\\phantom{-}1\\\phantom{-}1\end{bmatrix},
\begin{bmatrix}1\\3\\2\end{bmatrix},
\begin{bmatrix}1\\2\\1\end{bmatrix}
\right\}.
\end{equation*}
Note that $\begin{bmatrix}1 &3 & 2\end{bmatrix}^T \simeq k_1+k_2+k_4$,
and $\begin{bmatrix}1 &2 & 1\end{bmatrix}^T \simeq k_1+k_2+k_3$,
where $a\simeq b$ means that $a \, = \, \lambda \, b$ for some $\lambda > 0$.
Thus,~$A_iK \subset K, \, i=1,2$.  

To compute~$\check \rho(\cA)$ we note that~$\rho(A_1) = 1$ and hence, $\check \rho \le 1$. 
To prove that~$\check \rho \ge 1$, we construct a invariant {\em conic polytope}~$P \, = \, 
{co}\,(S) \, + \, K$, where $S$ is a finite subset of~$K\setminus\{0\}$. This set 
possesses the property~$A_i P \subset P, \, i=1,2$. 
Applying the  Invariant Polytope Algorithm~\cite{GP13} we get the invariant conic polytope~$P$
with  
\begin{align*}
S \, = \, \left\{
\begin{bmatrix}1\\1\\0\end{bmatrix},
\begin{bmatrix}-1\\\phantom{-}1\\\phantom{-}1\end{bmatrix},
\begin{bmatrix}2\\2\\2\end{bmatrix},
\begin{bmatrix}-2\\\phantom{-}4\\\phantom{-}2\end{bmatrix}
\right\}\, . 
\end{align*}
Let~$\bx$ be the closest to the origin point of~$P$. 
Then for every product~$\Pi = A(k)\cdots A(1)$, we have 
$\Pi\bx \in P$, hence,~$\|\Pi \bx\| \ge \|\bx\|$ and therefore, 
$\|\Pi\| \ge 1$. This implies that~$\check \rho \ge 1$. 

Thus,~$\check \rho(\cA) = 1$.  Several conclusions can be drawn from this equality. First of all, the 
linear switching system defined by the family~$\cA$ is not stabilisable. Second, 
the Lyapunov exponent~$\lambda(\cA)$, which is not less than~$\log \, \check \rho$, is nonnegative. 
Finally, the family~$\cA$ is not mortal. Note that all those results could hardly be obtained without solving the invariant cone problem for~$\cA$. We are not aware of any
algorithm of approximate computation of the lower spectral radius or of proving the non-mortality
for general set of matrices. 
We could solve those problems for matrices from Example~\ref{ex.10} because we constructed their common invariant cone.

}
\end{ex}

 \end{document}